\documentclass[11pt,reqno]{amsart}

\usepackage{cite,a4wide,amsmath,amssymb,graphicx,enumitem,amsaddr,tikz,pgfplots,hhline,epstopdf,multirow,multicol,soul}
\usepackage[hang,flushmargin]{footmisc}

\usepackage{tabu}

\usetikzlibrary{decorations.markings,patterns}

\newtheorem{theorem}{Theorem}

\setlist[enumerate]{leftmargin=8mm}
\setlist[itemize]{leftmargin=8mm}

\newcommand{\cR}{\mathcal{R}}

\newcommand{\e}{\text{e}}

\newcommand{\de}{\text{d}}

\definecolor{darkgreen}{RGB}{34,177,76}

\begin{document}
\title{Combating new COVID-19 variants in Indonesia: Is accelerating four-dose vaccinations sufficient?}

\author{\small Benny Yong, Jonathan Hoseana, Livia Owen}
\address{\normalfont\small Center for Mathematics and Society, Department of Mathematics, Parahyangan Catholic University, Bandung 40141, Indonesia}
\email{benny\_y@unpar.ac.id\textnormal{, }j.hoseana@unpar.ac.id\textnormal{, }livia.owen@unpar.ac.id}
\date{}

\begin{abstract}
As new COVID-19 variants continue to emerge globally, we develop an SISI-type mathematical model to determine whether the transmission risks arising when such a new variant enters Indonesia can be mitigated by solely accelerating the country's ongoing four-dose vaccination programme. We begin by determining the model's basic reproduction number, as well as the model's equilibria and their stability. Subsequently, employing parameter values representing the country's situation as of March 20, 2023, we conduct a numerical sensitivity analysis in two simulated cases corresponding to two different levels of the new variant's transmission. The results show that, a satisfactory mitigation relying solely on vaccinations necessitates a drastic acceleration in the low-transmission case, and proves unachievable in the high-transmission case. Accordingly, we recommend that the acceleration of the ongoing four-dose vaccinations be carried out in conjunction with other intervention measures, such as improvements of the vaccine's efficacy and the disease's recovery rate.


\end{abstract}

\maketitle

\section{Introduction}\label{section:Introduction}

Despite being declared no longer a global emergency by WHO in May 2023 \cite{WHO}, the Coronavirus Disease 2019 (COVID-19) remains a notable subject of public attention. In Indonesia, as in most countries worldwide, the disease has brought significant impact on both social and economic sectors \cite{OliviaGibsonNasrudin,SuryahadiIzzatiSuryadarma,SusilawatiFalefiPurwoko} as various virus variants successively entered the country. Indeed, the Alpha (B.1.1.7), Beta (B.1.351), and Delta (B.1.617.2) variants, first identified in the United Kingdom (September 2020), South Africa (May 2020), and India (October 2020), respectively, have all entered Indonesia in May 2021 \cite{MargariniAnindita}. Subsequently, the Omicron variant (B.1.1.529), which is reportedly more than four times as transmissible as Delta \cite{Matsumaya}, entered the country in December 2021 \cite{Wahyuni}. In October 2022, the Indonesian Ministry of Health confirmed the entry of two subvariants of Omicron \cite{Fernandez}: the XBB and BF.7 (BA.5.2.1.7) which, according to some reports, is the most transmissible and infectious subvariant of Omicron, possessing a shorter incubation period \cite{Mohammed,SinghAnandSrivastava}.

The COVID-19 pandemic has also brought a considerable scientific impact, with a surge of interest in studying the disease's spread, especially through the use of mathematical models. In 2022, we commenced our study by constructing a continuous SIR-type model which incorporates as a main parameter the susceptible individuals' cautiousness level, and using the model to justify mathematically that a high cautiousness level is crucial for the disease's eradication \cite{YongOwenHoseana}. 
Subsequently, employing a discretised version of the same model, and using the data of mid 2021, we studied the impact of social restrictions, developing a quantitative method to determine the optimal level of social restrictions to be implemented in Jakarta on any given day, based on the latest values of the disease's effective reproduction number and the hospitals' bed-occupancy rate \cite{YongHoseanaOwen1}.


More recently, we modified substantially our model in \cite{YongOwenHoseana} to take into account five forms of governmental interventions: vaccinations, social restrictions, tracings, testings, and treatments \cite{YongHoseanaOwen2}. The resulting model suggested that, in a disease-free state, when there is an increase in the number of new cases, the optimal intervention strategy is to reimplement social restrictions. On the other hand, in an endemic state, the optimal strategy for eradication is to administer vaccinations. However, effort should be directed not primarily towards the acceleration of vaccination rate, but towards the use of high-efficacy vaccines.

Since the initiation of Indonesia's national vaccination programme in January 2021 \cite{Triwadani}, the country has observed a tangible implementation of the above strategy, with over 70\% of the population having received their first two vaccine doses by September 2022 \cite{FirdausFardahHaryati}. In a response to such a progress, the Indonesian government decided to lift all COVID-19-related social restrictions at the end of 2022 \cite{Hidayat}. Such a decision sparked controversies, with concerns expressed by some epidemiologists regarding the potential risks involved. An argument supporting such concerns is that the currently low administration rate of boosters is not sufficient to overcome the risks associated with the aforementioned emergence of novel virus variants \cite{Mufarida,UlyaMovanita}.

This paper aims to determine an appropriate governmental response to this argument. More precisely, given that the Indonesian government has commenced the administration of fourth-dose vaccinations ---second booster shots--- in January 2023 \cite{MuthiarinyArkyasa} while novel variants continue to emerge worldwide \cite{Swaroop,Singh}, we aim to determine whether an acceleration of such four-dose COVID-19 vaccinations would be sufficient to mitigate the transmission risks which arise when a new variant enters Indonesia.

We shall begin by constructing an SISI-type mathematical model, which takes into account both fourth-dose vaccinations and infections by a new virus variant, and identifying a suitable positively invariant domain (section \ref{sec:model}). Subsequently, we determine the model's basic reproduction number, as well as the model's equilibria and their stability (section \ref{sec:equilibria}), and use the model to simulate two different cases corresponding to two different transmission levels of the new variant (section \ref{sec:numerical}). The results reveal that, in the lower transmission case, a significant intensification of four-dose vaccinations ---which must increase the percentage of vaccinated individuals by at least 13 times the current value--- is necessary for a successful mitigation, while in the higher transmission case, four-dose vaccinations alone ---even when intensified to target the entire population--- is insufficient for a success. Accordingly, we conclude that four-dose vaccinations must be administered in conjunction with other forms of intervention, such as increasing the vaccine's efficacy and the disease's recovery rate (section \ref{sec:conclusions}).

\section{The model}\label{sec:model}

Let us begin by constructing our SISI-type model. Once again, our aim is to determine the impact of four-dose vaccinations on the mitigation of risks associated with a new virus variant. Realising the suggestions proposed in \cite{OwenHoseanaYong,YongHoseanaOwen2}, our model shall feature reinfections \cite{OwenHoseanaYong} as well as a distinction between susceptible individuals who have received different numbers of vaccination doses \cite{YongHoseanaOwen2}. In addition, as a number of studies have confirmed the presence of antibody responses following a COVID-19 infection \cite{Altawalah,LiEtAl,SwartzEtAl,QiEtAl}, our model shall treat recovered and fully vaccinated susceptible individuals as having the same, higher level of immunity.

Accordingly, let us assume that the susceptible individuals can be divided into two different categories corresponding to two different immunity levels: a lower immunity level, possessed by those who have neither received their fourth-dose vaccination nor infected by the new variant, and a higher immunity level, possessed by the rest, i.e., susceptible individuals who have either received their fourth-dose vaccination or recovered from an infection of the new variant. On the other hand, let us also assume that the infected individuals can be divided into two different categories: those who are infected by the new variant for the first time, and those who are infected by the new variant but not for the first time. Our model, therefore, involves four different compartments $\text{S}$, $\text{I}_1$, $\text{S}_1$, $\text{I}_2$, and we denote by $S=S(t)$, $I_1=I_1(t)$, $S_1=S_1(t)$, and $I_2=I_2(t)$ the number of individuals in each of these compartments at time $t\geqslant0$ (Table \ref{tab:variables}).

\begin{table}\centering
\begin{tabular}{|c|l|}\hline
Variable & Description\\
\hhline{|=|=|}
$S=S(t)$ & \begin{minipage}{0.8\linewidth}\medskip Number of susceptible individuals who, at time $t$, have neither received their fourth-dose vaccination nor infected by the new variant\end{minipage}\\[0.4cm]\hline
$I_1=I_1(t)$ & \begin{minipage}{0.8\linewidth}\medskip Number of individuals who, at time $t$, are infected by the new variant for the first time\end{minipage}\\[0.4cm]\hline
$S_1=S_1(t)$ & \begin{minipage}{0.8\linewidth}\medskip Number of the remaining susceptible individuals, i.e., those who, at time $t$, have either received their fourth-dose vaccination or recovered from an infection of the new variant\end{minipage}\\[0.6cm]\hline
$I_2=I_2(t)$ & \begin{minipage}{0.8\linewidth}\medskip Number of individuals who, at time $t$, are infected by the new variant but not for the first time\end{minipage}\\[0.4cm]\hline
\end{tabular}
\caption{\label{tab:variables} Time-dependent variables involved in the model \eqref{eq:model}.}
\end{table}

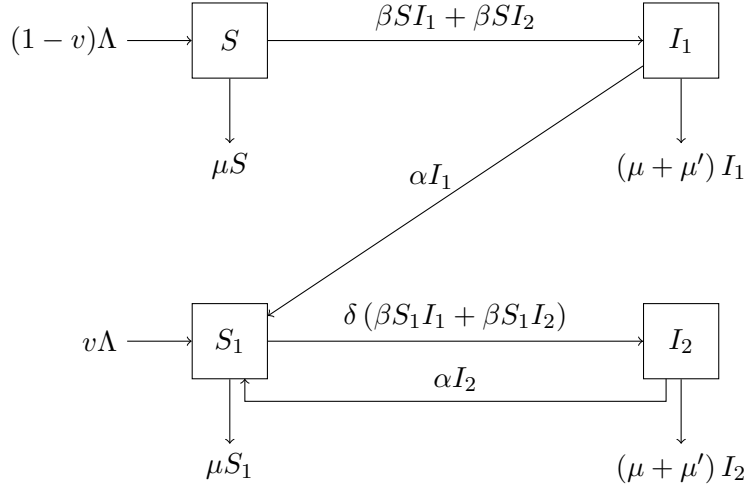
\begin{figure}\centering
\begin{tikzpicture}
\node[rectangle,draw=black,minimum size=1cm] (S) at (0,0) {$S$};
\node[rectangle,draw=black,minimum size=1cm] (I1) at (6,0) {$I_1$};
\node[rectangle,draw=black,minimum size=1cm] (S1) at (0,-4) {$S_1$};
\node[rectangle,draw=black,minimum size=1cm] (I2) at (6,-4) {$I_2$};

\draw[->] (S) -- (I1);
\node[above] at (3,0) {$\beta SI_1+\beta SI_2$};

\node (lahirS) at (-1.5,0) {};

\draw[->] (lahirS) -- (S);
\node[left=-4pt] at (lahirS) {$(1-v)\Lambda$};

\node (matiS) at (0,-1.5) {};

\draw[->] (S) -- (matiS);
\node[below=-4pt] at (matiS) {$\mu S$};

\node (matiI1) at (6,-1.5) {};

\draw[->] (I1) -- (matiI1);
\node[below=-4pt] at (matiI1) {$\left(\mu +\mu'\right)I_1$};

\draw[->] (I1) -- (S1);
\node[above left=-3pt] at (3,-2) {$\alpha I_1$};

\node (lahirS1) at (-1.5,-4) {};

\draw[->] (lahirS1) -- (S1);
\node[left=-4pt] at (lahirS1) {$v\Lambda$};

\draw[->] (S1) -- (I2);
\node[above] at (3,-4) {$\delta\left(\beta S_1I_1+\beta S_1I_2\right)$};

\node (matiS1) at (0,-5.5) {};

\draw[->] (S1) -- (matiS1);
\node[below=-4pt] at (matiS1) {$\mu S_1$};

\node (matiI2) at (6,-5.5) {};

\draw[->] (I2) -- (matiI2);
\node[below=-4pt] at (matiI2) {$\left(\mu +\mu'\right)I_2$};

\draw[->] (5.8,-4.5) -- (5.8,-4.8) -- (0.2,-4.8) -- (0.2,-4.5);

\node[above] at (3,-4.8) {$\alpha I_2$};
\end{tikzpicture}
\caption{\label{fig:compartement}The compartment diagram of our model \eqref{eq:model}.}
\end{figure}

To construct the model itself, let us first assume the recovery rates of individuals in $\text{I}_1$ and $\text{I}_2$ are linear: $\alpha_1 I$ and $\alpha_2 I_2$, respectively, where $\alpha_1,\alpha_2>0$, and that all  recovered individuals enter the compartment $\text{S}_1$. Next, we regard the infection rate of individuals in $\text{S}$ as being higher than that of individuals in $\text{S}_1$, as the latter is suppressed by the vaccine's efficacy. Assuming the infection rates to be bilinear, we thus formulate the transition rates from $\text{S}$ to $\text{I}_1$ and from $\text{S}_1$ to $\text{I}_2$ to be $\beta SI_1 + \beta SI_2$ and $\delta\left(\beta S_1I_1 + \beta S_1I_2\right)$, respectively, where $\beta>0$, $\delta\in[0,1]$ is such that $1-\delta$ represents the vaccine's efficacy, and $N=N(t)$ denotes the time-dependent total population:
\begin{equation}\label{eq:totalpopulation}
N=S+I_1+S_1+I_2.
\end{equation}
Finally, we assume that four-dose vaccinated individuals enter the population at the rate $v\Lambda$, while other individuals enter at the rate $(1-v)\Lambda$, where $\Lambda>0$ and $v\in[0,1]$ denote, respectively, the recruitment and vaccination rates. Employing the linear death rates $\mu S$, $\mu S_1$, $\left(\mu+\mu'\right)I_1$, and $\left(\mu+\mu'\right)I_2$ for the compartments $\text{S}$, $\text{S}_1$, $\text{I}_1$, and $\text{I}_2$, respectively, where $\mu,\mu'>0$, one constructs the compartment diagram in Figure \ref{fig:compartement}, and subsequently the model 
\begin{equation}\label{eq:model}
\left\{\begin{array}{rcl}
\dfrac{\de S}{\de t} \!\!\!&=&\!\!\!\left( 1-v \right) \Lambda -\mu S-\beta SI_{1}-\beta SI_{2},\\[0.3cm]
\dfrac{\de I_{1}}{\de t} \!\!\!&=&\!\!\!\beta SI_{1}+\beta SI_{2}-\mu
I_{1}-\mu ^{\prime }I_{1}-\alpha I_{1},\\[0.3cm]
\dfrac{\de S_{1}}{\de t} \!\!\!&=&\!\!\!v \Lambda +\alpha I_{1}+\alpha I_{2} -\mu S_{1}-\delta \left( 
\beta S_{1}I_{1}+\beta S_{1}I_{2}\right),\\[0.3cm]
\dfrac{\de I_{2}}{\de t} \!\!\!&=&\!\!\!\delta \left( \beta S_{1}I_{1}+\beta
S_{1}I_{2}\right)-\mu I_{2}-\mu ^{\prime }I_{2} -\alpha I_{2}.
\end{array}\right.
\end{equation}
For the reader's convenience, we present in Table \ref{tab:parameters} descriptions of the model's parameters and their values used in our numerical analysis. The values are chosen to represent the situation in Indonesia as of March 20, 2023.

\begin{table}\centering\renewcommand{\arraystretch}{2.5}
\begin{tabular}{|c|c|c|c|c|}
\hline
Parameter & Description & Unit & \begin{minipage}{2cm}\centering
Value for simulation\medskip
\end{minipage} &Source\\
\hhline{|=|=|=|=|=|}
$v $ & \begin{minipage}{3.8cm}\centering
proportion of vaccinated individuals\end{minipage} & dimensionless & \multicolumn{2}{c|}{varied; see section \ref{sec:numerical}}\\\hline
$\Lambda $ & recruitment rate & $\dfrac{\text{individual}}{\text{day}}$ & $%
\dfrac{234666020}{65\times 365}$ & \begin{minipage}{2.7cm}\smallskip\centering
estimated as $\mu N\left( 0\right)$,\\ $N\left( 0\right)$ denoting the initial target of vaccination \cite{AndiKurmalaNasution}\end{minipage} \\ \hline
$\mu $ & natural death rate & $\dfrac{1}{\text{day}}$ & $\dfrac{1}{65\times 365%
}$ & \cite{YongHoseanaOwen2}\\ \hline
$\beta $ & transmission rate & $\dfrac{1}{\text{individual}\times \text{day}%
}$ & \multicolumn{2}{c|}{varied; see section \ref{sec:numerical}}\\ \hline
$\mu ^{\prime }$ & \begin{minipage}{4cm}\centering
death rate due to COVID-19\end{minipage} & $\dfrac{1}{\text{day}}$ & $%
0.0291$ & \cite{YongHoseanaOwen2}\\ \hline
$\alpha $ & recovery rate & $\dfrac{1}{\text{day}}$ & $0.011$ & \cite{YongHoseanaOwen2}\\ \hline
$1-\delta $ & vaccine efficacy & dimensionless & $0.653$ & \cite{YongHoseanaOwen2}\\\hline
\end{tabular}\smallskip
\caption{\label{tab:parameters}Parameters involved in the model \eqref{eq:model} and their values used in our numerical analysis (section \ref{sec:numerical}).}
\end{table}

Now, let us consider a solution $\left(S(t),I_1(t),S_1(t),I_2(t)\right)$ of the model \eqref{eq:model} associated to an initial condition $\left(S(0),S_1(0),I_1(0),I_2(0)\right)\in\mathbb{R}_+^4$, where $\mathbb{R}_+:=[0,\infty)$. For every $t^\ast\geqslant 0$ satisfying $S\left(t^\ast\right)=0$, one obtains from the first equation of \eqref{eq:model} that
$$\left.\frac{\de S}{\de t}\right|_{t=t^\ast}=(1-v)\Lambda\geqslant 0,$$
which means that the function $S$ in non-decreasing at $t=t^\ast$. Since $S(0)\geqslant 0$, this implies that $S(t)\geqslant 0$ for every $t\geqslant 0$. Similarly, one proves that $I_1(t)\geqslant 0$, $S_1(t)\geqslant 0$, and $I_2(t)\geqslant 0$ for every $t\geqslant 0$.

Next, adding the four equations in the model \eqref{eq:model}, one finds that the total population \eqref{eq:totalpopulation} satisfies
$$\frac{\de N}{\de t}=\Lambda-\mu N -\mu' I_1 -\mu' I_2\leqslant \Lambda - \mu N,\qquad\text{i.e.,}\qquad \frac{\de N}{\de t}+\mu N\leqslant \Lambda.$$
Multiplying both sides $\e^{\mu t}$, one rewrites the last inequality as
$$\frac{\de }{\de t}\left[\e^{\mu t} N(t)\right]\leqslant \frac{\de }{\de t}\left[\frac{\Lambda}{\mu}\e^{\mu t}+N(0)-\frac{\Lambda}{\mu}\right].$$
Since $\e^{\mu t} N(t)$ and $(\Lambda/\mu)\e^{\mu t}+N(0)-\Lambda/\mu$ coincide at $t=0$, it follows that for every $t\geqslant 0$ we have
$$\e^{\mu t} N(t)\leqslant \frac{\Lambda}{\mu}\e^{\mu t}+N(0)-\frac{\Lambda}{\mu},\qquad\text{i.e.,}\qquad N(t)\leqslant \frac{\Lambda}{\mu}+\left[N(0)-\frac{\Lambda}{\mu}\right]\e^{-\mu t}\xrightarrow{t\to\infty}\frac{\Lambda}{\mu}.$$
We have therefore proved that the solution $\left(S(t),S_1(t),I_1(t),I_2(t)\right)$ is bounded, and established a bounded subset of $\mathbb{R}_+^4$ which is positively invariant under the model \eqref{eq:model}.\smallskip

\begin{theorem}\label{thm:domain}
The solution of the model \eqref{eq:model} is bounded for any initial condition in $\mathbb{R}_+^4$. Moreover, the set
$$\Omega:=\left\{\left(S,I_1,S_1,I_2\right)\in\mathbb{R}_+^4: S+I_1+S_1+I_2\leqslant\frac{\Lambda}{\mu}\right\}$$
is positively invariant under the model.
\end{theorem}\smallskip

\section{Equilibria and stability}\label{sec:equilibria}

Let us next study the equilibria of the model \eqref{eq:model} and their local stabilities. The equilibria satisfy the system
\begin{equation}\label{eq:equilibriasystem}
\left\{\begin{array}{rcl}
0 \!\!\!&=&\!\!\!f_1\left(S,I_1,S_1,I_2\right),\\[0.1cm]
0 \!\!\!&=&\!\!\!f_2\left(S,I_1,S_1,I_2\right),\\[0.1cm]
0 \!\!\!&=&\!\!\!f_3\left(S,I_1,S_1,I_2\right),\\[0.1cm]
0 \!\!\!&=&\!\!\!f_4\left(S,I_1,S_1,I_2\right),\end{array}\right.
\end{equation}
where
\begin{align*}
f_1\left(S,I_1,S_1,I_2\right)&=\left( 1-v \right) \Lambda -\mu S-\beta SI_{1}-\beta SI_{2},\\
f_2\left(S,I_1,S_1,I_2\right)&=\beta SI_{1}+\beta SI_{2}-\mu
I_{1}-\mu ^{\prime }I_{1}-\alpha I_{1},\\
f_3\left(S,I_1,S_1,I_2\right)&=v \Lambda +\alpha I_{1}+\alpha I_{2} -\mu S_{1}-\delta \left( 
\beta S_{1}I_{1}+\beta S_{1}I_{2}\right),\\
f_4\left(S,I_1,S_1,I_2\right)&=\delta \left( \beta S_{1}I_{1}+\beta
S_{1}I_{2}\right) -\alpha I_{2}-\mu I_{2}-\mu ^{\prime }I_{2}.
\end{align*}
Letting $I=I_1+I_2$, one obtains from each of the equations in \eqref{eq:equilibriasystem} that
\begin{equation}\label{eq:equilibria}
\left\{\begin{array}{rcl}
S\!\!\!&=&\!\!\!\dfrac{\left(1-\upsilon\right)\Lambda}{\mu+\beta I},\\[0.3cm]
I_1\!\!\!&=&\!\!\!\dfrac{\beta SI}{\mu+\mu^\prime+\alpha}=\dfrac{\beta\Lambda\left(1-\upsilon\right)I}{\left(\mu+\beta I\right)\left(\mu+\mu^\prime+\alpha\right)},\\[0.3cm]
S_1\!\!\!&=&\!\!\!\dfrac{\upsilon\Lambda + \alpha I}{\mu+\delta\beta I},\\[0.3cm]
I_2\!\!\!&=&\!\!\!\dfrac{\delta\beta S_1 I}{\mu+\mu^\prime+\alpha}=\dfrac{\delta\beta\left(\upsilon\Lambda+\alpha I\right)I}{\left(\mu+\delta\beta I\right)\left(\mu+\mu^\prime+\alpha\right)}.
\end{array}\right.
\end{equation}
Let us now consider two different cases: $I=0$ and $I>0$.

\subsection{The disease-free equilibrium}

Since $I_1,I_2\geqslant0$, in the case of $I=0$, i.e., $I_1+I_2=0$, we must have that $I_1=I_2=0$. Moreover, the first and third equations in \eqref{eq:equilibria} give
$$S=\frac{(1-v)\Lambda}{\mu}\qquad\text{and}\qquad S_1=\frac{v\Lambda}{\mu}.$$
We thus obtain the equilibrium
\begin{equation}\label{eq:DFE}
\mathbf{e}^{\left( 0\right) } =\left( S^{\left( 0\right) },I_{1}^{\left( 0\right)
},S_{1}^{\left( 0\right) },I_{2}^{\left( 0\right) }\right)=\left( \frac{\left( 1-v \right) \Lambda }{\mu },\,0,\,\frac{v
\Lambda }{\mu },\,0\right),
\end{equation}
i.e., the unique disease-free equilibrium of the model \eqref{eq:model}. To study its local stability, we notice that the Jacobian of the model \eqref{eq:model} is given by
\begin{align*}
\mathbf{J}\left(S,I_1,S_1,I_2\right) &=\left( 
\begin{array}{cccc}
\dfrac{\partial f_1}{\partial S} & \dfrac{\partial f_1}{\partial I_{1}} & \dfrac{\partial f_1}{\partial S_{1}} & \dfrac{\partial f_1}{\partial I_{2}} \\[0.3cm]
\dfrac{\partial f_2}{\partial S} & \dfrac{\partial f_2
}{\partial I_{1}} & \dfrac{\partial f_2}{\partial S_{1}} & \dfrac{f_2}{\partial I_{2}} \\ [0.3cm]
\dfrac{\partial f_3}{\partial S} & \dfrac{\partial f_3}{\partial I_{1}} & \dfrac{\partial f_3}{\partial S_{1}} & \dfrac{\partial f_3}{\partial I_{2}} \\[0.3cm] 
\dfrac{\partial f_4}{\partial S} & \dfrac{\partial f_4}{\partial I_{1}} & \dfrac{\partial f_4}{\partial S_{1}} & \dfrac{\partial f_4}{\partial I_{2}}
\end{array}\right) \\[0.1cm]
&=\renewcommand{\arraycolsep}{1.5pt}\left( 
\begin{array}{cccc}
-\mu -\beta I_{1}-\beta I_{2} & -\beta S & 0 & -\beta S \\[0.15cm]
\beta I_{1}+\beta I_{2} & \beta S-\mu -\mu ^{\prime }-\alpha & 0 & \beta S
\\[0.15cm] 
0 & \alpha -\delta \beta S_{1} & -\mu -\delta \left( \beta I_{1}+\beta
I_{2}\right) & \alpha-\delta \beta S_{1} \\ [0.15cm]
0 & \delta \beta S_{1} & \delta \left( \beta I_{1}+\beta I_{2}\right) & 
\delta \beta S_{1} -\mu -\mu ^{\prime }-\alpha
\end{array}%
\right).
\end{align*}
Evaluating the Jacobian at the disease-free equilibrium $\mathbf{e}^{\left(0\right)}$ gives the matrix
\begin{align*}
\mathbf{J}\left(\mathbf{e}^{\left(0\right)}\right) &=\left( 
\begin{array}{cccc}
-\mu & - \dfrac{\Lambda\beta\left( 1-\upsilon \right) }{\mu } & 0 & -\dfrac{\Lambda\beta\left( 1-\upsilon \right) }{\mu } \\[0.3cm]
0 & \dfrac{\Lambda\beta\left( 1-\upsilon \right) }{\mu }-\mu -\mu ^{\prime
}-\alpha & 0 & \dfrac{\Lambda\beta\left( 1-\upsilon \right) }{\mu } \\[0.3cm]
0 & \alpha -\dfrac{\Lambda \beta\delta \upsilon  }{\mu } & -\mu & \alpha-\dfrac{\Lambda \beta\delta \upsilon  }{\mu } \\[0.3cm]
0 & \dfrac{\Lambda \beta\delta \upsilon  }{\mu } & 0 & \dfrac{\Lambda \beta\delta \upsilon  }{\mu } -\mu -\mu ^{\prime }-\alpha
\end{array}
\right)
\end{align*}
whose eigenvalues are
$$\lambda _{1}=\lambda_2=-\mu,\quad\lambda _{3}= -\mu -\mu ^{\prime }-\alpha,\quad \text{and}\quad \lambda _{4}=\frac{\Lambda\beta \left(\delta \upsilon+1-\upsilon\right)}{\mu }-\mu-\mu^{\prime}-\alpha.$$
Since $\lambda _i<0$ for all $i\in\{1,2,3\}$, then the disease-free equilibrium $\mathbf{e}^{\left(0\right)}$ is locally asymptotically stable if $\lambda_4<0$, i.e., if $\cR_0<1$, and is unstable if $\lambda_4<0$, i.e., if $\cR_{0}>1$, where
\begin{equation}\label{eq:R0}
\cR_0:=\frac{\Lambda \beta \left( \delta
\upsilon +1-\upsilon \right) }{\mu \left( \mu +\mu ^{\prime }+\alpha\right) 
}
\end{equation}
denotes the model's \textit{basic reproduction number}, obtainable via, e.g., the next-generation approach \cite{DiekmannHeesterbeekMetz,DriesscheWatmough}. We have therefore proved the following theorem.\smallskip

\begin{theorem}\label{thm:DFE}
For all sets of parameter values, the model \eqref{eq:model} has a unique disease-free equilibrium $\mathbf{e}^{\left(0\right)}$ given by \eqref{eq:DFE}. This equilibrium is locally asymptotically stable if $\cR_0<1$, and is unstable if $\cR_{0}>1$.
\end{theorem}\smallskip


\subsection{The endemic equilibria}

Let us next consider the case $I>0$. In this case, the equation $I=I_1 + I_2$, i.e.,
$$I=\dfrac{\beta\Lambda\left(1-\upsilon\right)I}{\left(\mu+\beta I\right)\left(\mu+\mu^\prime+\alpha\right)}+\dfrac{\delta\beta\left(\upsilon\Lambda+\alpha I\right)I}{\left(\mu+\delta\beta I\right)\left(\mu+\mu^\prime+\alpha\right)}$$
is equivalent to the quadratic equation
\begin{equation}\label{eq:QE}
\mathcal{A}I^2 + \mathcal{B}I + \mathcal{C}=0,
\end{equation}
where
\begin{align}
\mathcal{A} &=\delta \beta ^{2}\left( \mu +\mu ^{\prime }\right), \label{eq:A} \\
\mathcal{B} &=\beta \left[ \delta \mu \left( \mu +\mu ^{\prime }\right) +\mu\left(\mu +\mu ^{\prime }+\alpha\right) -\delta
\beta \Lambda \right]\nonumber\\
&=\frac{\beta \left[ \beta \Lambda \left[ \left( 1-\delta \right) \left(
1-\upsilon \right) +\delta \left( 1-\cR_0\right)
\right] +\cR_0\mu \delta \left( \mu +\mu ^{\prime
}\right) \right] }{\cR_0}, \label{eq:B}\\
\mathcal{C} &=\mu ^{2}\left(\mu +\mu ^{\prime }+\alpha\right)\left( 1-\cR_0\right).\label{eq:C}
\end{align}
Now, if $\mathcal{R}_0>1$, then $\mathcal{A}>0$ and $\mathcal{C}<0$, implying the existence of a unique positive endemic equilibrium. On the other hand, if $\mathcal{R}_0<1$, then $\mathcal{A}>0$, $\mathcal{B}>0$, and $\mathcal{C}>0$, implying that no positive endemic equilibria exist. Finally, if $\mathcal{R}_0=1$, then $\mathcal{C}=0$, and so \eqref{eq:QE} is equivalent to $I\left(\mathcal{A} I +\mathcal{B}\right)=0$, which means that $I=0$ or $I=-\mathcal{B}/\mathcal{A}<0$, and thus no positive endemic equilibria exist. In summary, we have the following theorem.\smallskip

\begin{theorem}\label{thm:EE}
The model \eqref{eq:model} has a unique positive endemic equilibrium $$\mathbf{e}^{\left(1\right)}=\left( S^{\left( 1\right) },I_{1}^{\left( 1\right)
},S_{1}^{\left( 1\right) },I_{2}^{\left( 1\right) }\right),$$ 
where
\begin{eqnarray*}
&S^{(1)}=\dfrac{\left(1-\upsilon\right)\Lambda}{\mu+\beta I^{(1)}},\qquad I_1^{(1)}=\dfrac{\beta\Lambda\left(1-\upsilon\right)I^{(1)}}{\left[\mu+\beta I^{(1)}\right]\left(\mu+\mu^\prime+\alpha\right)},\qquad S_1^{(1)}=\dfrac{\upsilon\Lambda + \alpha I^{(1)}}{\mu+\delta\beta I^{(1)}},&\\
&I_2^{(1)}=\dfrac{\delta\beta\left[\upsilon\Lambda+\alpha I^{(1)}\right]I^{(1)}}{\left[\mu+\delta\beta I^{(1)}\right]\left(\mu+\mu^\prime+\alpha\right)},&
\end{eqnarray*}
and $I^{(1)}=I_1^{(1)}+I_2^{(1)}>0$, if and only if $\mathcal{R}_0>1$.
\end{theorem}\smallskip

Building upon Theorem \ref{thm:EE}, it is natural to hypothesise that if $\mathcal{R}_0>1$, then the positive endemic equilibrium $\mathbf{e}^{\left(1\right)}$, which exist uniquely, is locally asymptotically stable. However, proving this analytically requires tedious symbolic computation of the model's Jacobian at the endemic equilibrium $\mathbf{e}^{\left(1\right)}$, whose explicit expression, obtainable from \eqref{eq:QE} via the quadratic formula, is already complicated. At this point, therefore, let us migrate from analytical to numerical methods, with which we shall confirm graphically the above hypothesis for our specified set of parameter values, and subsequently accomplish our initial goal, i.e., to determine whether the acceleration of four-dose vaccinations is sufficient for countering new variant transmissions.

\section{Numerical analysis}\label{sec:numerical}

As mentioned in section \ref{sec:model}, we shall use for our numerical analysis the values of the parameters $\Lambda$, $\mu$, $\mu'$, $\alpha$, and $1-\delta$ which are listed in Table \ref{tab:parameters}. Once again, these values are specified to represent the situation in Indonesia as of March 20, 2023.

We shall present our numerical analysis in two subsections. The first subsection comprises a preliminary discussion on the model's basic reproduction number, the local asymptotic stability of the model's endemic equilibrium, and the transcritical bifurcation occurring at the basic reproduction threshold. The second subsection comprises a two-stage sensitivity analysis, of the basic reproduction number at the first stage and of the endemic-equilibrium subpopulation sizes at the second stage, carried out with the aim of answering our main research question, i.e., whether accelerating four-dose vaccinations alone suffices to overcome the new variant's transmission risks.

\subsection{Basic reproduction number analysis}
\begin{figure}\centering
\begin{tikzpicture}
\pgfset{declare function={f(\x)=(2226980529800000/380949)*10^(-10)*\x*(-653/1000*(167/10000)+1);}}

\begin{axis}[
	xmin=0,
	xmax=4.5,
	ymin=0,
	ymax=2.602,
	xtick={0,1.73},
	xticklabels={{\small $0$},{\small $\beta_{0}$}},
	ytick={0.5,1,1.5,2,2.5},
	yticklabels={{\small $0.50$},{\small $1.00$},{\small $1.50$},{\small $2.00$},{\small $2.50$}},
	axis on top=true,
	samples=100,
	xlabel=$\beta$,
	ylabel=$\cR_0(\beta)$,
	width=6.75cm,
	height=6.75cm,
	scaled x ticks=false,
	ylabel near ticks
]
\addplot [very thick,domain=0:4.5] {f(x)};

\draw[dashed] (axis cs:1.73,0) -- (axis cs:1.73,1);
\draw[dashed] (axis cs:1.73,1) -- (axis cs:0,1);
\end{axis}
\end{tikzpicture}
\caption{\label{fig:R0}The graph of $\mathcal{R}_0\left(\beta\right)$ versus $\beta$ for $\upsilon=\upsilon_1=1.67\%$.}
\end{figure}
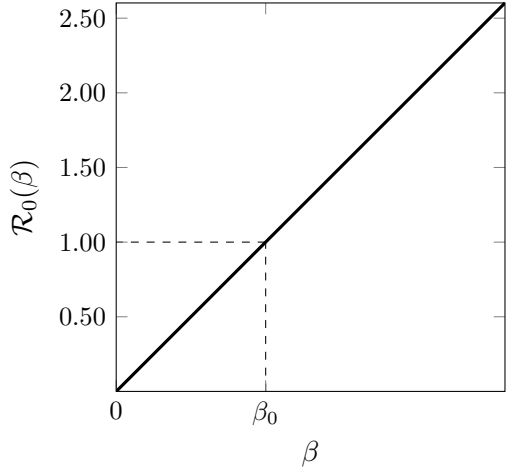

For the aforementioned parameter values, the model's basic reproduction number \eqref{eq:R0} is given by
\begin{equation}\label{eq:R0dlmbetaupsilon}
\mathcal{R}_0\left(\beta,\upsilon\right)=\frac{2226980529800000}{380949}\beta\left(1-\frac{653}{1000}\upsilon\right).
\end{equation}
A plot of the curve $\mathcal{R}_0\left(\beta,\upsilon\right)=1$ is presented in Figure \ref{fig:R0=1}. The shaded region is the feasible disease-free region, i.e., the region comprising pairs $(\beta,\upsilon)\in[0,\infty)\times [0,1]$ for which $\cR_0(\beta,\upsilon)<1$. One immediately sees that administering four-dose vaccinations to at least $\upsilon_1=25\%$, $\upsilon_2=50\%$, and $\upsilon_3=75\%$ of the entering individuals, respectively, suffices to overcome the new variant's transmissions occurring at the rates of $\beta=\beta_1$, $\beta=\beta_2$, and $\beta=\beta_3$, where
$$\beta_1\approx 2.04\cdot 10^{-10},\qquad \beta_2\approx 2.54\cdot 10^{-10},\qquad\text{and}\qquad \beta_3\approx 3.35\cdot 10^{-10}.$$
In reality, however, as of March 20, 2023, four-dose vaccinations have only been administered to a very small fraction of Indonesia's population, namely, $\upsilon=\upsilon_0=0.0167=1.67\%$ \cite{SatuanTugas}, for which we have that
\begin{equation}\label{eq:R0dlmbeta}
\mathcal{R}_0\left(\beta\right)=\mathcal{R}_0\left(\beta,\upsilon_0\right)=\frac{36711584740407967}{6349150}\beta,
\end{equation}
as plotted in Figure \ref{fig:R0}. We see that $\cR_0<1$ and $\cR_0>1$ correspond respectively to $\beta<\beta_0$ and $\beta>\beta_0$, where $$\beta_0\approx 1.73\cdot 10^{-10}.$$
That is, the proportion of vaccinated individuals achieved on March 20, 2023 suffices only to mitigate transmissions occurring at the rate not exceeding $1.73\cdot 10^{-10}$ per individual per day.

\begin{figure}\centering
\begin{tikzpicture}
\pgfset{declare function={f(\x)=(1/1454218285959400)*(-380949+2226980529800000*\x*10^(-10))/(\x*10^(-10));}}

\begin{axis}[
	xmin=0,
	xmax=4.9,
	ymin=0,
	ymax=1,
	xtick={0,2.04,2.54,3.35},
	xticklabels={{\small $0$},{\small $\beta_{1}$},{\small $\beta_{2}$},{\small $\beta_{3}$}},
	ytick={0.25,0.5,0.75,1},
	yticklabels={{\small $0.25$},{\small $0.50$},{\small $0.75$},{\small $1.00$}},
	axis on top=true,
	samples=100,
	xlabel=$\beta$,
	ylabel=$\upsilon$,
	width=6.75cm,
	height=6.75cm,
	scaled x ticks=false,
	ylabel near ticks
]
\fill[pattern=north east lines] (axis cs:0,0) rectangle (axis cs:4.9,1);

\addplot+[mark=none,draw=none,fill=white,domain=0:4.9] {f(x)} \closedcycle;
\addplot [very thick,domain=0:4.9] {f(x)};
\node at (axis cs:3.25,0.495) {$\mathcal{R}_0=1$};

\draw[dashed] (axis cs:2.04,0) -- (axis cs:2.04,0.25);
\draw[dashed] (axis cs:2.54,0) -- (axis cs:2.54,0.5);
\draw[dashed] (axis cs:3.35,0) -- (axis cs:3.35,0.75);

\draw[dashed] (axis cs:0,0.25) -- (axis cs:2.04,0.25);
\draw[dashed] (axis cs:0,0.5) -- (axis cs:2.54,0.5);
\draw[dashed] (axis cs:0,0.75) -- (axis cs:3.35,0.75);
\end{axis}
\end{tikzpicture}
\caption{\label{fig:R0=1}The curve $\cR_0=1$ on the $\beta\upsilon$-plane, with the disease-free region shaded.}
\end{figure}

\begin{figure}
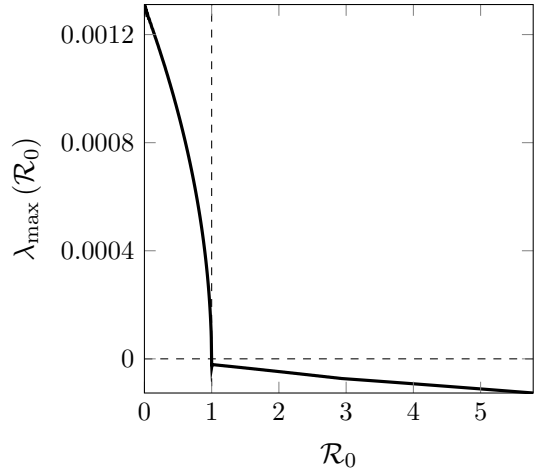
\centering

\caption{\label{fig:lambdamax}Plot of $\lambda_{\max}\left(\cR_0\right)$ versus $\cR_0$ for $\cR_0\in[0.95,1.05]$, indicating that $\mathbf{e}^{(1)}\left(\cR_0\right)$ is locally asymptotically stable if $\cR_0>1$, and is unstable if $\cR_0<1$.}
\end{figure}

Let us now confirm the local asymptotic stability of $\mathbf{e}^{(1)}$ in the case of $\cR_0>1$, i.e., $\beta>\beta_0$. From \eqref{eq:R0dlmbeta} one obtains that
\begin{equation}\label{eq:betadlmR0}
\beta\left(\mathcal{R}_0\right)=\frac{6349150}{36711584740407967}\cR_0.
\end{equation}
Substituting this, $\upsilon=\upsilon_0$, and the other parameters' values in Table \ref{tab:parameters} into \eqref{eq:A}, \eqref{eq:B}, and \eqref{eq:C} gives the following exact expressions for the coefficients of \eqref{eq:QE} in terms of $\cR_0$:
\begin{align*}
\mathcal{A}\left(\cR_0\right)&=\frac{1547420825179752197}{5116022763961580751185838014489445844000}{\cR_0}^2,\\
\mathcal{B}\left(\cR_0\right)&=\frac{11190529508566}{109006467232912135109843300161}\cR_0\left(\frac{1572380558426459}{440631010000000}-\cR_0\right),\\
\mathcal{C}\left(\cR_0\right)&=\frac{380949}{5341689681250000}\left(1-\cR_0\right).
\end{align*}
The quadratic formula subsequently provides an explicit expression for the total number of infected individuals at the model's unique positive endemic equilibrium $\mathbf{e}^{(1)}\left(\cR_0\right)$:
\begin{equation}\label{eq:endemicI}
I^{(1)}\left(\cR_0\right)=\frac{-\mathcal{B}\left(\cR_0\right)+\sqrt{\left[\mathcal{B}\left(\cR_0\right)\right]^2-4\mathcal{A}\left(\cR_0\right)\mathcal{C}\left(\cR_0\right)}}{2\mathcal{A}\left(\cR_0\right)},
\end{equation}
and consequently those of the equilibrium's coordinates $S^{(1)}\left(\cR_0\right)$, $I_1^{(1)}\left(\cR_0\right)$, $S_1^{(1)}\left(\cR_0\right)$, and $I_2^{(1)}\left(\cR_0\right)$, via the formulae provided by Theorem \ref{thm:EE}. Using these to compute the model's Jacobian at $\mathbf{e}^{(1)}\left(\cR_0\right)$,
$$\mathbf{J}\left(\mathbf{e}^{(1)}\left(\cR_0\right)\right)= \mathbf{J}\left(S^{(1)}\left(\cR_0\right),\,I_1^{(1)}\left(\cR_0\right),\,S_1^{(1)}\left(\cR_0\right),\,I_2^{(1)}\left(\cR_0\right)\right),$$
one observes graphically (Figure \ref{fig:lambdamax}) that if $\cR_0>1$ then $\lambda_{\max}\left(\cR_0\right)<0$, where $\lambda_{\max}\left(\cR_0\right)$ denotes the maximum real part of the eigenvalues of $\mathbf{J}\left(\mathbf{e}^{(1)}\left(\cR_0\right)\right)$, thereby confirming that $\mathbf{e}^{(1)}\left(\cR_0\right)$ is locally asymptotically stable if $\cR_0>1$.

\begin{figure}
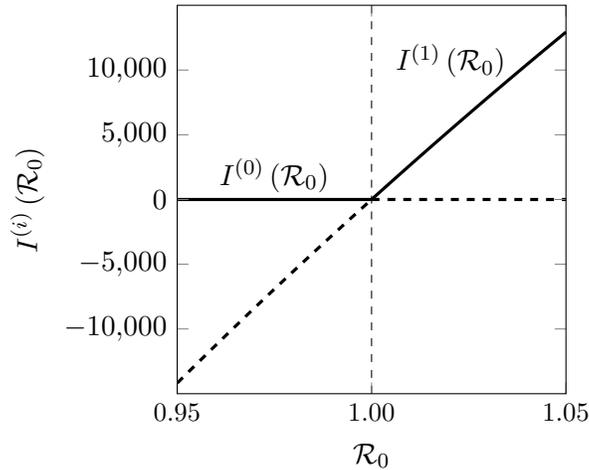
\centering
%
\caption{\label{fig:transcritical}Plot of $I^{(i)}\left(\cR_0\right)$ versus $\cR_0$ for $\cR_0\in[0.95,1.05]$. Solid and dashed lines for $I^{(i)}\left(\cR_0\right)$ indicate, respectively, local asymptotic stability and instability of the equilibrium $\mathbf{e}^{(i)}\left(\cR_0\right)$. Notice the transcritical bifurcation occurring at $\cR_0=1$.}
\end{figure}

Furthermore, in the case of $\cR_0<1$, we observe the following. First, since $\mathcal{A}\left(\cR_0\right)>0$, $\mathcal{B}\left(\cR_0\right)>0$, and $\mathcal{C}\left(\cR_0\right)>0$, then $I^{(1)}\left(\cR_0\right)<0$. Second, since $\lambda_{\max}\left(\cR_0\right)>0$ (Figure \ref{fig:lambdamax}), then $\mathbf{e}^{(1)}\left(\cR_0\right)$ is unstable. A plot of $I^{(1)}\left(\cR_0\right)$ and $I^{(0)}\left(\cR_0\right)=0$, versus $\cR_0$, is shown in Figure \ref{fig:transcritical}, where a solid (dashed) line for $I^{(i)}\left(\cR_0\right)$ indicates local asymptotic stability (instability) of $\mathbf{e}^{(i)}\left(\cR_0\right)$. As we can see, the model \eqref{eq:model} undergoes a transcritical bifurcation, i.e., a stability exchange between the equilibria $\mathbf{e}^{(0)}\left(\cR_0\right)$ and $\mathbf{e}^{(1)}\left(\cR_0\right)$, at $\cR_0=1$.

\subsection{Sensitivity analysis}

In the previous subsection we have found that, at the state where four-dose COVID-19 vaccinations are administered only to $\upsilon=\upsilon_0=1.67\%$ of the Indonesian population, the risks posed by a new variant are surmountable only if transmissions occur at a sufficiently low rate, $\beta<\beta_0$, where $\beta_0\approx 1.73\cdot 10^{-10}$. In this subsection, we shall focus on the endemic case $\beta>\beta_0$, and carry out a two-stage sensitivity analysis \cite{ChitnisHymanCushing} to quantitatively assess the sensitivity of, firstly, the model's basic reproduction number, and secondly, the subpopulation sizes at the positive endemic equilibrium, with respect to small changes of each of the model's parameters.

\sloppy For the first stage of our analysis, let us recall that the sensitivity of the basic reproduction number $\cR_{0}$ to small changes of a parameter $p\in\left\{\upsilon,\Lambda,\mu,\beta,\mu',\alpha,\delta\right\}$ is quantified by the sensitivity index
$$\Upsilon _{p}^{\cR_0}:=\frac{\partial \cR
_{0}}{\partial p}\times \frac{p}{\cR
_{0}}$$
of $\cR_{0}$ with respect to $p$, thus defined assuming the required differentiability \cite[page 1280]{ChitnisHymanCushing}. Straightforward computation leads to the expressions
\begin{eqnarray*}
&\Upsilon _{\upsilon }^{\cR_0} =\dfrac{(\delta -1)\upsilon }{(\delta
-1)\upsilon +1},\qquad
\Upsilon _{\Lambda }^{\cR_0} =1,\qquad
\Upsilon _{\mu }^{\cR_0} =-\dfrac{\mu }{\alpha +\mu +\mu ^{\prime }}-1,\qquad\Upsilon _{\beta }^{\cR_0} =1,&\\
&\Upsilon _{\mu ^{\prime }}^{\cR_0} = \dfrac{-\mu ^{\prime }}{\alpha +\mu +\mu ^{\prime }},\qquad\Upsilon _{\alpha }^{\cR_0} =\dfrac{-\alpha }{\alpha +\mu +\mu ^{\prime }},\qquad\text{and}\qquad\Upsilon _{\delta }^{\cR_0} =\dfrac{\delta \upsilon }{\delta \upsilon
+1-\upsilon }.&
\end{eqnarray*}
Notice that none of these sensitivity indices depend on $\beta$. Substituting the aforementioned values of parameters (i.e., those in Table \ref{tab:parameters} and $\upsilon=0.0167$) gives the values presented in Table \ref{tab:sensitivityparameters}.

\begin{table}\renewcommand{\arraystretch}{1.2}
\begin{tabular}{|c|r|r|}\hline
$p$ & \multicolumn{1}{c|}{$\Upsilon _{p}^{\cR_0}$} & \multicolumn{1}{c|}{$1/\Upsilon _{p}^{\cR_0}$} \\\hhline{|=|=|=|}
$\upsilon $ & $-0.01103$ & $+90.70021$ \\ \hline
$\Lambda $ & $+1.00000$ & $-1.00000$ \\ \hline
$\mu $ & $-1.00105$ & $+0.99895$ \\ \hline
$\beta $ & $+1.00000$ & $-1.00000$ \\ \hline
$\mu ^{\prime }$ & $-0.72492$ & $+1.37946$ \\ \hline
$\alpha $ & $-0.27403$ & $+3.64929$ \\ \hline
$\delta $ & $+0.00586$ & $-170.68369$\\\hline
\end{tabular}\smallskip
\caption{\label{tab:sensitivityparameters}The sensitivity indices of $\cR_0$ with respect to $p\in\left\{\upsilon,\Lambda,\mu,\beta,\mu',\alpha,\delta\right\}$ for parameter values given by Table \ref{tab:parameters} and $\upsilon=0.0167$. None of these indices depend on $\beta$.}
\end{table}

Therefore, in estimation, increasing the proportion $\upsilon$ of vaccinated individuals by 10\% suppresses $\cR_0$ merely by $0.1103\%$, while increasing the vaccine efficacy $1-\delta $ by 10\% suppresses $\cR_0$ merely by $0.0586\%$. Indeed, to achieve a mere 1\% suppression of $\cR_0$, it is necessary to increase either the proportion of vaccinated individuals alone by $90.7002\%$, or the vaccine efficacy alone by $170.6837\%$, both of which appearing impractical. This shows that neither the proportion of vaccinated individuals nor the vaccine efficacy hold sole significance in controlling transmission risks. In fact, as apparent from Table \ref{tab:sensitivityparameters}, these are the two parameters upon which the number $\cR_0$ depends least significantly. The number $\cR_0$ depends most significantly, besides on the recruitment and natural death rates $\Lambda$ and $\mu$ which are practically unchangeable, on the transmission rate $\beta$, confirming that reductions of contacts between susceptible and infected individuals, i.e., social restrictions, remain an effective form of intervention. However, in the present situation whereby the economy has started to recover after almost two years of social restrictions \cite{SatyagrahaKenzu,SuroyoSulaiman,SwasonoSuhanda}, reinforcement of social restrictions might not be desirable.

To determine a suitable course of action which relies on four-dose vaccinations, let us proceed with the second stage of our sensitivity analysis. Let $\mathcal{P}:=\left\{\upsilon,\Lambda,\mu,\beta,\mu',\alpha,\delta\right\}$. Fix a parameter $p\in\mathcal{P}$. We aim to assess the sensitivity of the subpopulation sizes $S^{(1)}$, $I_1^{(1)}$, $S_1^{(1)}$, and $I_2^{(1)}$ at the model's endemic equilibrium $\mathbf{e}^{(1)}$, with respect to $p$. As noted in \cite[page 1294]{ChitnisHymanCushing}, this is quantified by the sensitivity indices
\begin{eqnarray*}
&\Upsilon^{S^{(1)}}_{p}=\left.\dfrac{\partial S}{\partial p}\right|_{S=S^{(1)}}\cdot\dfrac{p}{S^{(1)}},\qquad \Upsilon^{I_1^{(1)}}_{p}=\left.\dfrac{\partial I_1}{\partial p}\right|_{I_1=I_1^{(1)}}\cdot\dfrac{p}{I_1^{(1)}},\qquad \Upsilon^{S_1^{(1)}}_{p}=\left.\dfrac{\partial S_1}{\partial p}\right|_{S_1=S_1^{(1)}}\cdot\dfrac{p}{S_1^{(1)}},&\\
&\text{and}\qquad \Upsilon^{I_2^{(1)}}_{p}=\left.\dfrac{\partial I_2}{\partial p}\right|_{I_2=I_2^{(1)}}\cdot\dfrac{p}{I_2^{(1)}}&
\end{eqnarray*}
of $S^{(1)}$, $I_1^{(1)}$, $S_1^{(1)}$, and $I_2^{(1)}$ with respect to $p$.

To compute the above indices, we first need to compute 
\begin{equation}\label{eq:partial}
\left.\dfrac{\partial S}{\partial p}\right|_{S=S^{(1)}},\quad \left.\dfrac{\partial I_1}{\partial p}\right|_{I_1=I_1^{(1)}},\quad \left.\dfrac{\partial S_1}{\partial p}\right|_{S_1=S_1^{(1)}}\left.\dfrac{\partial I_2}{\partial p}\right|_{I_2=I_2^{(1)}},\quad\text{and}\quad \left.\dfrac{\partial I_2}{\partial p}\right|_{I_2=I_2^{(1)}}.
\end{equation}
For this purpose, for every $i\in\{1,2,3,4\}$, we differentiate both sides of the equation $0=f_i=f_i\left(S,I_1,S_1,I_2\right)$ with respect to $p$, obtaining, by the multivariate chain rule,
$$0=\frac{\partial f_i}{\partial p}=\frac{\partial f_i}{\partial S}\cdot\frac{\partial S}{\partial p}+\frac{\partial f_i}{\partial I_1}\cdot\frac{\partial I_1}{\partial p}+\frac{\partial f_i}{\partial S_1}\cdot\frac{\partial S_1}{\partial p}+\frac{\partial f_i}{\partial I_2}\cdot\frac{\partial I_2}{\partial p}+\sum_{p'\in\mathcal{P}}\frac{\partial f_i}{\partial p'}\cdot\frac{\partial p'}{\partial p},$$
which, since
$$\frac{\partial p'}{\partial p}=\begin{cases}
0,&\text{if }p'\neq p;\\
1,&\text{if }p'=p,
\end{cases}$$
is equivalent to
$$\frac{\partial f_i}{\partial S}\cdot\frac{\partial S}{\partial p}+\frac{\partial f_i}{\partial I_1}\cdot\frac{\partial I_1}{\partial p}+\frac{\partial f_i}{\partial S_1}\cdot\frac{\partial S_1}{\partial p}+\frac{\partial f_i}{\partial I_2}\cdot\frac{\partial I_2}{\partial p}=-\frac{\partial f_i}{\partial p}.$$
Since $i\in\{1,2,3,4\}$, we have a system of four equations, which can be written in a matrix form as
$$\mathbf{J}\left(S,I_1,S_1,S_2\right)\mathbf{X}_p=\mathbf{Y}_p,$$
where
$$\mathbf{X}_p:=\left(\begin{array}{c}
 \dfrac{\partial S}{\partial p}\\[0.3cm]
 \dfrac{\partial I_1}{\partial p}\\[0.3cm]
 \dfrac{\partial S_1}{\partial p}\\[0.3cm]
 \dfrac{\partial I_2}{\partial p}
 \end{array}\right)\qquad\text{and}\qquad
 \mathbf{Y}_p:=\left(\begin{array}{c}
 \dfrac{\partial f_1}{\partial p}\\[0.3cm]
 \dfrac{\partial f_2}{\partial p}\\[0.3cm]
 \dfrac{\partial f_3}{\partial p}\\[0.3cm]
 \dfrac{\partial f_4}{\partial p}
 \end{array}\right).$$
Thus, the desired quantities \eqref{eq:partial} can be obtained by evaluating the vector
$$\mathbf{X}_p=\left[\mathbf{J}\left(S,I_1,S_1,I_2\right)\right]^{-1}\mathbf{Y}_p$$
at $\left(S,I_1,S_1,I_2\right)=\mathbf{e}^{(1)}$, after which the four sensitivity indices $\Upsilon^{S^{(1)}}_{p}$, $\Upsilon^{I_1^{(1)}}_{p}$, $\Upsilon^{S_1^{(1)}}_{p}$, and $\Upsilon^{I_2^{(1)}}_{p}$ are readily computable.

Given the unavailability of the data of the actual transmission rate $\beta$, let us carry out the above computation of $\Upsilon^{S^{(1)}}_{p}$, $\Upsilon^{I_1^{(1)}}_{p}$, $\Upsilon^{S_1^{(1)}}_{p}$, and $\Upsilon^{I_2^{(1)}}_{p}$ for each $p\in\mathcal{P}$ in two simulated cases, corresponding to two qualitatively different transmission levels \cite[sec.\ 3]{ChitnisHymanCushing}: a low-transmission case $\beta=\beta_l=2\cdot 10^{-10}$ and a high-transmission case $\beta=\beta_h=8\cdot 10^{-10}$.\smallskip

\subsubsection*{The low-transmission case}

For the parameter values presented in Table \ref{tab:parameters}, $\upsilon=0.0167$, and $\beta=2\cdot 10^{-10}$, we have that $\cR_0\approx 1.15643$ and
$$\mathbf{e}^{(1)}=\left(S^{(1)},I_1^{(1)},S_1^{(1)}, I_2^{(1)}\right)\approx\left(1.96\cdot 10^8,\, 3.62\cdot 10^4,\, 1.28\cdot 10^7,\, 8.19\cdot 10^{2}\right).$$

\begin{figure}
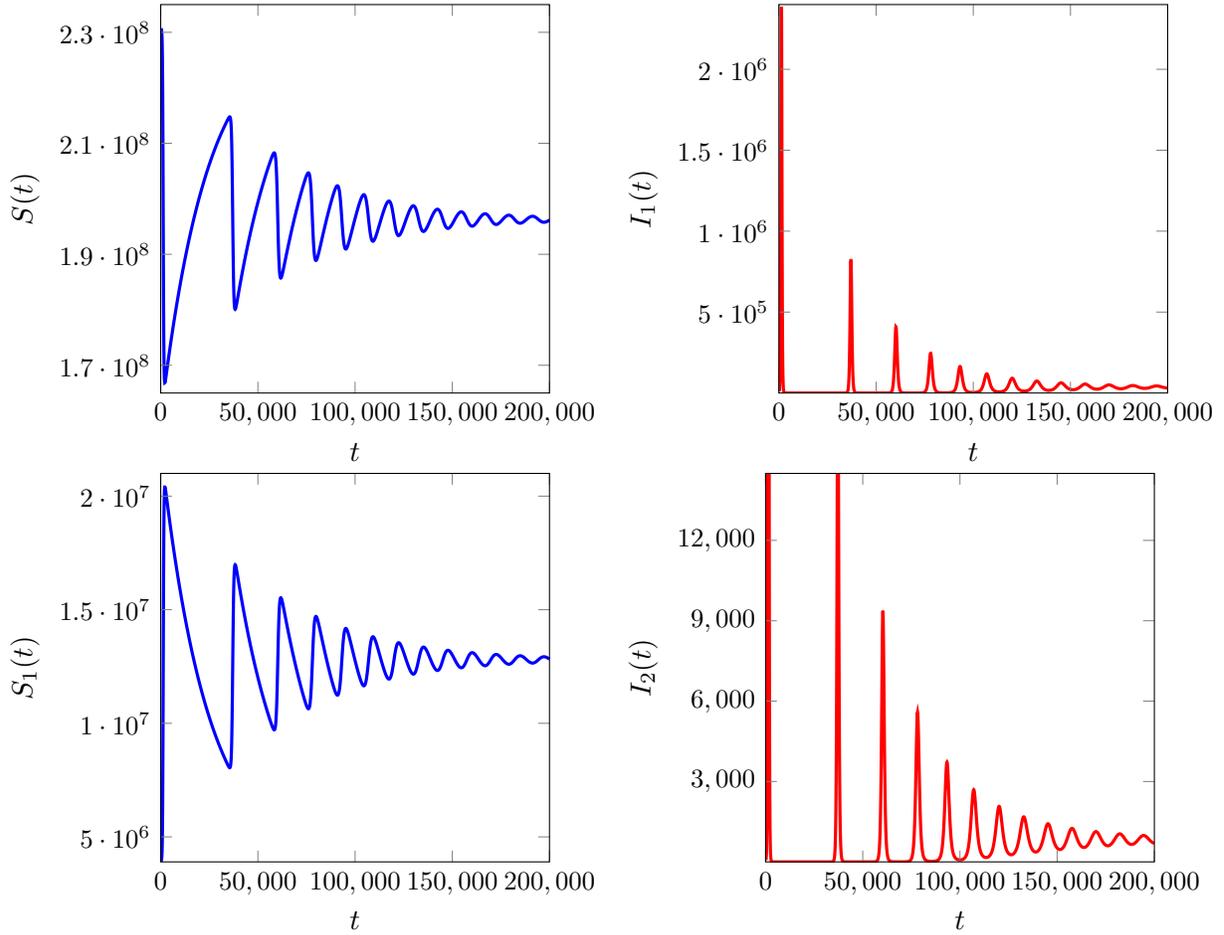
\centering
%
\caption{\label{fig:TSlow} Plots of $S(t)$, $I_1(t)$, $S_1(t)$, and $I_2(t)$, versus $t$, for $0\leqslant t\leqslant 200,000$, in the low-transmission case.}
\end{figure}

Let us specify an initial condition representing the situation in Indonesia on March 20, 2023:
\begin{equation}\label{eq:IC}
\left\{
\begin{array}{rcl}
 S\left( 0\right) \!\!\!&=&\!\!\!(1-\upsilon) N\left( 0\right)-I_{1}\left( 0\right)=98.33\%\cdot 234666020-3660=230743437,\\[0.1cm]
I_{1}\left( 0\right)\!\!\!&=&\!\!\!3660,\\[0.1cm]
S_{1}\left( 0\right) \!\!\!&=&\!\!\!\upsilon N\left( 0\right)= 1.67\%\cdot 234666020=3918922,\\[0.1cm]
 I_2(0)\!\!\!&=&\!\!\!0,
 \end{array}\right.
\end{equation}
where, as in Table \ref{tab:parameters}, $N(0)$ denotes the initial target of vaccination \cite{AndiKurmalaNasution,Pinandhita}. Figure \ref{fig:TSlow} shows the convergence of the model's solution associated to the initial condition given by \eqref{eq:IC} to the endemic equilibrium $\mathbf{e}^{(1)}$. Notice that convergence is achieved via decreasing-amplitude oscillations.

\begin{table}\renewcommand{\arraystretch}{1.2}
\begin{tabular}{|c|r|r|r|r|}\hline
$p$ & \multicolumn{1}{c|}{$\Upsilon _{p}^{S^{(1)}}$} & \multicolumn{1}{c|}{$\Upsilon _{p}^{I_1^{(1)}}$} & \multicolumn{1}{c|}{$\Upsilon _{p}^{S_1^{(1)}}$} & \multicolumn{1}{c|}{$\Upsilon _{p}^{I_2^{(1)}}$}\\ \hhline{|=|=|=|=|=|}
$\upsilon $ &  $-0.00538$      &$ -0.08303 $      & $0.23779$         & $0.16014$ \\\hline
$\mu $ & $0.11838$         & $-6.36751$         & $-5.18350$         & $-11.66939$            \\\hline
$\beta $ &$ -1.11719$        & $6.35968$ & $4.17828$          & $11.65515$          \\\hline
$\mu'$ & $0.82283$         & $-5.40895$       &$ -3.60137$         & $-9.83315$        \\\hline
$\alpha $ & $0.29317$          & $-1.94289$          &$-0.57169$        & $-2.80775$         \\\hline
$\delta $  & $-0.02368$       & $0.13478$         & $0.04617$        & $1.20462$ \\\hline
$\Lambda $ & $-0.11719$        & $7.35968$          & $5.17828$ & $12.65515$   \\\hline
\end{tabular}\smallskip
\caption{\label{tab:sensitivitylow}The sensitivity indices of $S^{(1)}$, $I_1^{(1)}$, $S_1^{(1)}$, and $I_2^{(1)}$ with respect to $p\in\left\{\upsilon,\Lambda,\mu,\beta,\mu',\alpha,\delta\right\}$ in the low-transmission case.}
\end{table}

\begin{figure}\centering
	\includegraphics[width=0.5\textwidth]{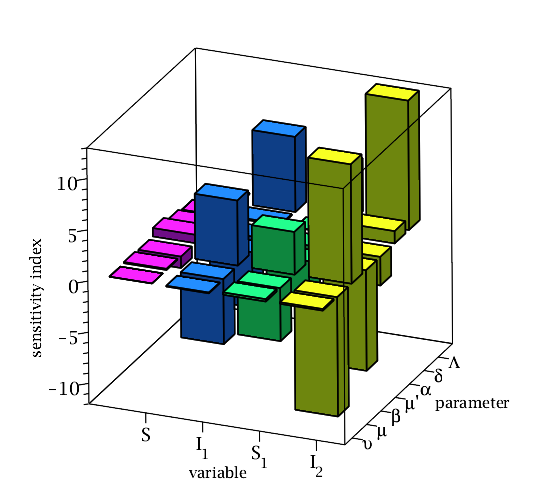}
	\caption{\label{fig:sensitivitylow}Three-dimensional boxplot of the sensitivity indices of $S^{(1)}$, $I_1^{(1)}$, $S_1^{(1)}$, and $I_2^{(1)}$ with respect to $p\in\left\{\upsilon,\Lambda,\mu,\beta,\mu',\alpha,\delta\right\}$ in the low-transmission case.}
\end{figure}

In this case, the values of the sensitivity indices are presented in Table \ref{tab:sensitivitylow} and visualised in Figure \ref{fig:sensitivitylow}. Notice that $I_1^{(1)}$ and $I_2^{(1)}$ are not only exceptionally sensitive to the birth, death, and transmission rates $\Lambda$, $\mu$, and $\beta$, but also insensitive to the proportion $\upsilon$ of vaccinated individuals and to the vaccine efficacy $1-\delta$, confirming our previous observation (Table \ref{tab:sensitivityparameters}). This means that, at the state of equilibrium, achieving a significant reduction of the number of infected individuals by relying only on vaccinations, albeit is not impossible, requires a considerable effort. Indeed, while four-dose vaccinations have been administered only to $\upsilon=1.67\%$ of the Indonesian population as of March 20, 2023, in this low-transmission case, one finds from \eqref{eq:R0dlmbetaupsilon} that it is necessary to administer four-dose vaccinations to at least $22.16\%$ of the population in order to achieve the disease-free state $\cR_0<1$. That is, the administration of four-dose vaccinations must be intensified until the proportion of fully vaccinated individuals surpasses 13 times the current percentage.

\subsubsection*{The high-transmission case}

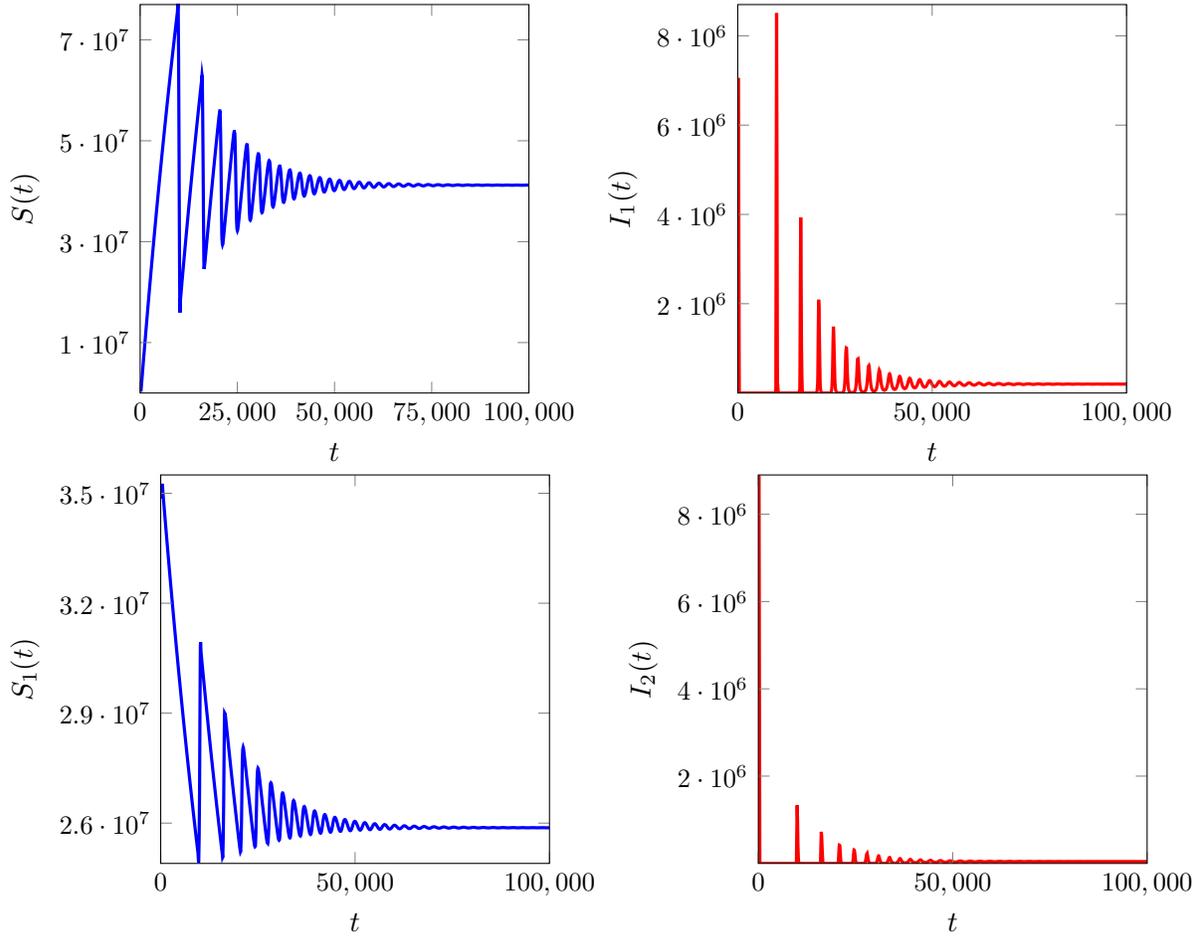
\begin{figure}\centering
\begin{tikzpicture}
\begin{axis}[
	xmin=0,
	xmax=100000,
	ymin=0,
	ymax=7.7,
	xtick={0,25000,50000,75000,100000},
	xticklabels={{\small $0$},{\small $25,000$},{\small $50,000$},{\small $75,000$},{\small $100,000$}},
	ytick={1,3,5,7},
	yticklabels={{\small $1\cdot 10^7$},{\small $3\cdot 10^7$},{\small $5\cdot 10^7$},{\small $7\cdot 10^7$}},
	axis on top=true,
	samples=300,
	xlabel=$t$,
	ylabel=$S(t)$,
	width=6.75cm,
	height=6.75cm,
	scaled x ticks=false,
	ylabel near ticks
]
\addplot [very thick,smooth,blue] coordinates {(190., 0.25415e-1) (390., .20106) (590., .39302) (790., .58346) (990., .77231) (1190., .95956) (1390., 1.1452) (1590., 1.3294) (1790., 1.5120) (1990., 1.6930) (2190., 1.8725) (2390., 2.0505) (2590., 2.2271) (2790., 2.4021) (2990., 2.5757) (3190., 2.7478) (3390., 2.9185) (3590., 3.0877) (3790., 3.2555) (3990., 3.4219) (4190., 3.5869) (4390., 3.7506) (4590., 3.9128) (4790., 4.0737) (4990., 4.2332) (5190., 4.3914) (5390., 4.5483) (5590., 4.7039) (5790., 4.8581) (5990., 5.0111) (6190., 5.1627) (6390., 5.3131) (6590., 5.4622) (6790., 5.6101) (6990., 5.7568) (7190., 5.9022) (7390., 6.0464) (7590., 6.1893) (7790., 6.3311) (7990., 6.4717) (8190., 6.6111) (8390., 6.7493) (8590., 6.8864) (8790., 7.0223) (8990., 7.1571) (9190., 7.2907) (9390., 7.4233) (9590., 7.5545) (9790., 7.6550) (9990., 5.1124) (10190., 1.7813) (10390., 1.8695) (10590., 2.0452) (10790., 2.2217) (10990., 2.3968) (11190., 2.5704) (11390., 2.7426) (11590., 2.9133) (11790., 3.0826) (11990., 3.2504) (12190., 3.4169) (12390., 3.5820) (12590., 3.7456) (12790., 3.9079) (12990., 4.0688) (13190., 4.2284) (13390., 4.3866) (13590., 4.5436) (13790., 4.6992) (13990., 4.8534) (14190., 5.0064) (14390., 5.1581) (14590., 5.3086) (14790., 5.4577) (14990., 5.6056) (15190., 5.7523) (15390., 5.8978) (15590., 6.0416) (15790., 6.1788) (15990., 6.1964) (16190., 4.5893) (16390., 2.6104) (16590., 2.5442) (16790., 2.6969) (16990., 2.8659) (17190., 3.0353) (17390., 3.2035) (17590., 3.3703) (17790., 3.5358) (17990., 3.6998) (18190., 3.8625) (18390., 4.0238) (18590., 4.1838) (18790., 4.3424) (18990., 4.4997) (19190., 4.6556) (19390., 4.8103) (19590., 4.9636) (19790., 5.1156) (19990., 5.2659) (20190., 5.4129) (20390., 5.5424) (20590., 5.5299) (20790., 4.6321) (20990., 3.2015) (21190., 2.9385) (21390., 3.0383) (21590., 3.1935) (21790., 3.3573) (21990., 3.5219) (22190., 3.6856) (22390., 3.8482) (22590., 4.0095) (22790., 4.1694) (22990., 4.3280) (23190., 4.4852) (23390., 4.6408) (23590., 4.7945) (23790., 4.9447) (23990., 5.0851) (24190., 5.1829) (24390., 5.0728) (24590., 4.3191) (24790., 3.4377) (24990., 3.2206) (25190., 3.2904) (25390., 3.4293) (25590., 3.5846) (25790., 3.7439) (25990., 3.9037) (26190., 4.0627) (26390., 4.2205) (26590., 4.3764) (26790., 4.5298) (26990., 4.6783) (27190., 4.8150) (27390., 4.9151) (27590., 4.8891) (27790., 4.5160) (27990., 3.8416) (28190., 3.4682) (28390., 3.4318) (28590., 3.5260) (28790., 3.6612) (28990., 3.8095) (29190., 3.9618) (29390., 4.1145) (29590., 4.2656) (29790., 4.4128) (29990., 4.5520) (30190., 4.6722) (30390., 4.7424) (30590., 4.6802) (30790., 4.3618) (30990., 3.8927) (31190., 3.6178) (31390., 3.5764) (31590., 3.6499) (31790., 3.7683) (31990., 3.9035) (32190., 4.0446) (32390., 4.1859) (32590., 4.3231) (32790., 4.4503) (32990., 4.5546) (33190., 4.6069) (33390., 4.5466) (33590., 4.3087) (33790., 3.9715) (33990., 3.7447) (34190., 3.6874) (34390., 3.7354) (34590., 3.8330) (34790., 3.9518) (34990., 4.0782) (35190., 4.2046) (35390., 4.3240) (35590., 4.4264) (35790., 4.4930) (35990., 4.4890) (36190., 4.3697) (36390., 4.1402) (36590., 3.9124) (36790., 3.7931) (36990., 3.7838) (37190., 3.8420) (37390., 3.9344) (37590., 4.0419) (37790., 4.1531) (37990., 4.2594) (38190., 4.3504) (38390., 4.4106) (38590., 4.4150) (38790., 4.3355) (38990., 4.1736) (39190., 3.9945) (39390., 3.8784) (39590., 3.8480) (39790., 3.8817) (39990., 3.9528) (40190., 4.0422) (40390., 4.1374) (40590., 4.2285) (40790., 4.3048) (40990., 4.3520) (41190., 4.3511) (41390., 4.2853) (41590., 4.1617) (41790., 4.0262) (41990., 3.9332) (42190., 3.9035) (42390., 3.9266) (42590., 3.9835) (42790., 4.0580) (42990., 4.1382) (43190., 4.2139) (43390., 4.2745) (43590., 4.3071) (43790., 4.2976) (43990., 4.2380) (44190., 4.1398) (44390., 4.0381) (44590., 3.9693) (44790., 3.9478) (44990., 3.9677) (45190., 4.0154) (45390., 4.0780) (45590., 4.1447) (45790., 4.2057) (45990., 4.2512) (46190., 4.2701) (46390., 4.2531) (46590., 4.1984) (46790., 4.1198) (46990., 4.0442) (47190., 3.9957) (47390., 3.9834) (47590., 4.0028) (47790., 4.0441) (47990., 4.0967) (48190., 4.1513) (48390., 4.1990) (48590., 4.2309) (48790., 4.2385) (48990., 4.2165) (49190., 4.1674) (49390., 4.1051) (49590., 4.0497) (49790., 4.0171) (49990., 4.0124) (50190., 4.0320) (50390., 4.0680) (50590., 4.1121) (50790., 4.1560) (50990., 4.1919) (51190., 4.2123) (51390., 4.2114) (51590., 4.1871) (51790., 4.1444) (51990., 4.0959) (52190., 4.0561) (52390., 4.0354) (52590., 4.0364) (52790., 4.0557) (52990., 4.0871) (53190., 4.1235) (53390., 4.1580) (53590., 4.1838) (53790., 4.1952) (53990., 4.1885) (54190., 4.1642) (54390., 4.1283) (54590., 4.0912) (54790., 4.0633) (54990., 4.0513) (55190., 4.0561) (55390., 4.0747) (55590., 4.1017) (55790., 4.1313) (55990., 4.1575) (56190., 4.1750) (56390., 4.1797) (56590., 4.1697) (56790., 4.1470) (56990., 4.1176) (57190., 4.0898) (57390., 4.0710) (57590., 4.0651) (57790., 4.0722) (57990., 4.0895) (58190., 4.1123) (58390., 4.1358) (58590., 4.1552) (58790., 4.1662) (58990., 4.1662) (59190., 4.1547) (59390., 4.1345) (59590., 4.1111) (59790., 4.0907) (59990., 4.0786) (60190., 4.0770) (60390., 4.0852) (60590., 4.1007) (60790., 4.1197) (60990., 4.1380) (61190., 4.1517) (61390., 4.1578) (61590., 4.1547) (61790., 4.1431) (61990., 4.1258) (62190., 4.1076) (62390., 4.0931) (62590., 4.0858) (62790., 4.0869) (62990., 4.0955) (63190., 4.1091) (63390., 4.1245) (63590., 4.1383) (63790., 4.1475) (63990., 4.1501) (64190., 4.1453) (64390., 4.1343) (64590., 4.1200) (64790., 4.1061) (64990., 4.0962) (65190., 4.0924) (65390., 4.0951) (65590., 4.1034) (65790., 4.1150) (65990., 4.1273) (66190., 4.1374) (66390., 4.1432) (66590., 4.1433) (66790., 4.1378) (66990., 4.1280) (67190., 4.1164) (67390., 4.1061) (67590., 4.0996) (67790., 4.0981) (67990., 4.1018) (68190., 4.1094) (68390., 4.1191) (68590., 4.1286) (68790., 4.1358) (68990., 4.1390) (69190., 4.1376) (69390., 4.1319) (69590., 4.1235) (69790., 4.1144) (69990., 4.1070) (70190., 4.1030) (70390., 4.1031) (70590., 4.1071) (70790., 4.1139) (70990., 4.1218) (71190., 4.1289) (71390., 4.1337) (71590., 4.1352) (71790., 4.1329) (71990., 4.1275) (72190., 4.1205) (72390., 4.1135) (72590., 4.1083) (72790., 4.1061) (72990., 4.1072) (73190., 4.1112) (73390., 4.1171) (73590., 4.1233) (73790., 4.1286) (73990., 4.1316) (74190., 4.1318) (74390., 4.1291) (74590., 4.1243) (74790., 4.1186) (74990., 4.1133) (75190., 4.1099) (75390., 4.1089) (75590., 4.1106) (75790., 4.1144) (75990., 4.1192) (76190., 4.1241) (76390., 4.1278) (76590., 4.1295) (76790., 4.1289) (76990., 4.1262) (77190., 4.1220) (77390., 4.1174) (77590., 4.1136) (77790., 4.1115) (77990., 4.1114) (78190., 4.1133) (78390., 4.1167) (78590., 4.1206) (78790., 4.1243) (78990., 4.1268) (79190., 4.1276) (79390., 4.1266) (79590., 4.1239) (79790., 4.1204) (79990., 4.1169) (80190., 4.1142) (80390., 4.1130) (80590., 4.1135) (80790., 4.1154) (80990., 4.1183) (81190., 4.1215) (81390., 4.1242) (81590., 4.1258) (81790., 4.1260) (81990., 4.1247) (82190., 4.1223) (82390., 4.1194) (82590., 4.1168) (82790., 4.1150) (82990., 4.1144) (83190., 4.1152) (83390., 4.1170) (83590., 4.1194) (83790., 4.1219) (83990., 4.1238) (84190., 4.1248) (84390., 4.1245) (84590., 4.1232) (84790., 4.1211) (84990., 4.1188) (85190., 4.1169) (85390., 4.1157) (85590., 4.1156) (85790., 4.1165) (85990., 4.1182) (86190., 4.1202) (86390., 4.1221) (86590., 4.1234) (86790., 4.1238) (86990., 4.1233) (87190., 4.1221) (87390., 4.1203) (87590., 4.1185) (87790., 4.1172) (87990., 4.1165) (88190., 4.1166) (88390., 4.1176) (88590., 4.1190) (88790., 4.1206) (88990., 4.1220) (89190., 4.1229) (89390., 4.1230) (89590., 4.1224) (89790., 4.1212) (89990., 4.1198) (90190., 4.1184) (90390., 4.1175) (90590., 4.1172) (90790., 4.1175) (90990., 4.1184) (91190., 4.1196) (91390., 4.1209) (91590., 4.1219) (91790., 4.1224) (91990., 4.1223) (92190., 4.1216) (92390., 4.1206) (92590., 4.1195) (92790., 4.1185) (92990., 4.1179) (93190., 4.1178) (93390., 4.1182) (93590., 4.1190) (93790., 4.1200) (93990., 4.1210) (94190., 4.1217) (94390., 4.1219) (94590., 4.1217) (94790., 4.1211) (94990., 4.1202) (95190., 4.1193) (95390., 4.1186) (95590., 4.1182) (95790., 4.1183) (95990., 4.1187) (96190., 4.1194) (96390., 4.1203) (96590., 4.1210) (96790., 4.1214) (96990., 4.1215) (97190., 4.1212) (97390., 4.1206) (97590., 4.1199) (97790., 4.1192) (97990., 4.1187) (98190., 4.1186) (98390., 4.1187) (98590., 4.1192) (98790., 4.1198) (98990., 4.1204) (99190., 4.1209) (99390., 4.1212) (99590., 4.1212) (99790., 4.1208) (99990., 4.1203)};
\end{axis}
\end{tikzpicture}\,\,\,\begin{tikzpicture}
\begin{axis}[
	xmin=0,
	xmax=100000,
	ymin=0,
	ymax=8.7,
	xtick={0,50000,100000,150000,200000},
	xticklabels={{\small $0$},{\small $50,000$},{\small $100,000$},{\small $150,000$},{\small $200,000$}},
	ytick={2,4,6,8},
	yticklabels={{\small $2\cdot 10^6$},{\small $4\cdot 10^6$},{\small $6\cdot 10^6$},{\small $8\cdot 10^6$}},
	axis on top=true,
	samples=300,
	xlabel=$t$,
	ylabel=$I_1(t)$,
	width=6.75cm,
	height=6.75cm,
	scaled x ticks=false,
	ylabel near ticks
]
\addplot [very thick,smooth,red] coordinates {(190., 7.05881) (390., 0) (590., 0.437011e-5) (790., 0) (990., 0) (1190., 0) (1390., 0) (1590., 0) (1790., 0) (1990., 0) (2190., 0) (2390., 0) (2590., 0) (2790., 0) (2990., 0) (3190., 0) (3390., 0) (3590., 0) (3790., 0) (3990., 0) (4190., 0) (4390., 0) (4590., 0) (4790., 0) (4990., 0) (5190., 0) (5390., 0) (5590., 0) (5790., 0) (5990., 0) (6190., 0) (6390., 0) (6590., 0) (6790., 0) (6990., 0) (7190., 0) (7390., 0) (7590., 0) (7790., 0) (7990., 0) (8190., 0) (8390., 0) (8590., 0) (8790., 0) (8990., 0) (9190., 0) (9390., 0) (9590., 0) (9790., .123254) (9990., 8.50667) (10190., .727860) (10390., 0.164828e-1) (10590., 0) (10790., 0) (10990., 0) (11190., 0) (11390., 0) (11590., 0) (11790., 0) (11990., 0) (12190., 0) (12390., 0) (12590., 0) (12790., 0) (12990., 0) (13190., 0) (13390., 0) (13590., 0) (13790., 0) (13990., 0) (14190., 0) (14390., 0) (14590., 0) (14790., 0) (14990., 0) (15190., 0) (15390., 0) (15590., 0) (15790., 0.178076e-1) (15990., .381965) (16190., 3.93115) (16390., 1.09725) (16590., 0.861398e-1) (16790., 0) (16990., 0) (17190., 0) (17390., 0) (17590., 0) (17790., 0) (17990., 0) (18190., 0) (18390., 0) (18590., 0) (18790., 0) (18990., 0) (19190., 0) (19390., 0) (19590., 0) (19790., 0) (19990., 0) (20190., 0) (20390., 0.463015e-1) (20590., .399576) (20790., 2.08994) (20990., 1.33406) (21190., .229641) (21390., 0.386744e-1) (21590., 0) (21790., 0) (21990., 0) (22190., 0) (22390., 0) (22590., 0) (22790., 0) (22990., 0) (23190., 0) (23390., 0) (23590., 0) (23790., 0) (23990., 0.243208e-1) (24190., .114277) (24390., .551996) (24590., 1.48728) (24790., .951833) (24990., .258978) (25190., 0.667561e-1) (25390., 0.207599e-1) (25590., 0) (25790., 0) (25990., 0) (26190., 0) (26390., 0) (26590., 0) (26790., 0) (26990., 0.121637e-1) (27190., 0.321336e-1) (27390., .101340) (27590., .343177) (27790., .893424) (27990., .989243) (28190., .454860) (28390., .156949) (28590., 0.582257e-1) (28790., 0.261689e-1) (28990., 0.147588e-1) (29190., 0.105411e-1) (29390., 0) (29590., 0.108869e-1) (29790., 0.155809e-1) (29990., 0.276492e-1) (30190., 0.595765e-1) (30390., .148484) (30590., .377319) (30790., .729272) (30990., .740843) (31190., .408083) (31390., .178507) (31590., 0.815638e-1) (31790., 0.437191e-1) (31990., 0.286835e-1) (32190., 0.233576e-1) (32390., 0.236506e-1) (32590., 0.296220e-1) (32790., 0.453247e-1) (32990., 0.826092e-1) (33190., .169921) (33390., .350362) (33590., .578906) (33790., .601064) (33990., .393877) (34190., .208101) (34390., .110557) (34590., 0.663448e-1) (34790., 0.473053e-1) (34990., 0.408537e-1) (35190., 0.428823e-1) (35390., 0.543162e-1) (35590., 0.814432e-1) (35790., .138973) (35990., .249619) (36190., .410109) (36390., .512536) (36590., .441895) (36790., .289079) (36990., .172211) (37190., .107412) (37390., 0.757053e-1) (37590., 0.624612e-1) (37790., 0.610712e-1) (37990., 0.705821e-1) (38190., 0.948670e-1) (38390., .143198) (38590., .227543) (38790., .342229) (38990., .426287) (39190., .403733) (39390., .302047) (39590., .202089) (39790., .136301) (39990., .100176) (40190., 0.836008e-1) (40390., 0.805629e-1) (40590., 0.896380e-1) (40790., .113415) (40990., .157812) (41190., .227807) (41390., .312750) (41590., .370111) (41790., .356060) (41990., .285579) (42190., .208410) (42390., .151861) (42590., .118145) (42790., .101994) (42990., 0.993501e-1) (43190., .109176) (43390., .133250) (43590., .174770) (43790., .233511) (43990., .295835) (44190., .330962) (44390., .315748) (44590., .263170) (44790., .204630) (44990., .159311) (45190., .131028) (45390., .117578) (45590., .116706) (45790., .127890) (45990., .152159) (46190., .190290) (46390., .238390) (46590., .282397) (46790., .301082) (46990., .284001) (47190., .242926) (47390., .198338) (47590., .163121) (47790., .140873) (47990., .131003) (48190., .132454) (48390., .144979) (48590., .168764) (48790., .202554) (48990., .240385) (49190., .269779) (49390., .277235) (49590., .259555) (49790., .226796) (49990., .192644) (50190., .165755) (50390., .149095) (50590., .142824) (50790., .146524) (50990., .159917) (51190., .182319) (51390., .211015) (51590., .239423) (51790., .257686) (51990., .258043) (52190., .241004) (52390., .214716) (52590., .188507) (52790., .168297) (52990., .156427) (53190., .153336) (53390., .158837) (53590., .172410) (53790., .192614) (53990., .215944) (54190., .236254) (54390., .246476) (54590., .242726) (54790., .227126) (54990., .206057) (55190., .186003) (55390., .171099) (55590., .163159) (55790., .162640) (55990., .169329) (56190., .182390) (56390., .199853) (56590., .218032) (56790., .231799) (56990., .236507) (57190., .230672) (57390., .216905) (57590., .200113) (57790., .184889) (57990., .174179) (58190., .169352) (58390., .170746) (58590., .178000) (58790., .189988) (58990., .204492) (59190., .218085) (59390., .226837) (59590., .227961) (59790., .221334) (59990., .209515) (60190., .196249) (60390., .184842) (60590., .177417) (60790., .174969) (60990., .177653) (61190., .184928) (61390., .195479) (61590., .207087) (61790., .216834) (61990., .221917) (62190., .220853) (62390., .214220) (62590., .204296) (62790., .193936) (62990., .185548) (63190., .180663) (63390., .179955) (63590., .183390) (63790., .190272) (63990., .199216) (64190., .208186) (64390., .214861) (64590., .217372) (64790., .215090) (64990., .208896) (65190., .200724) (65390., .192751) (65590., .186739) (65790., .183777) (65990., .184273) (66190., .188024) (66390., .194239) (66590., .201570) (66790., .208266) (66990., .212585) (67190., .213371) (67390., .210520) (67590., .204992) (67790., .198386) (67990., .192357) (68190., .188196) (68390., .186651) (68590., .187917) (68790., .191662) (68990., .197061) (69190., .202884) (69390., .207705) (69590., .210283) (69790., .209975) (69990., .206972) (70190., .202200) (70390., .196956) (70590., .192497) (70790., .189754) (70990., .189217) (71190., .190914) (71390., .194431) (71590., .198967) (71790., .203456) (71990., .206787) (72190., .208125) (72390., .207173) (72590., .204275) (72790., .200265) (72990., .196182) (73190., .192971) (73390., .191293) (73590., .191439) (73790., .193314) (73990., .196470) (74190., .200169) (74390., .203525) (74590., .205712) (74790., .206199) (74990., .204918) (75190., .202271) (75390., .198981) (75590., .195867) (75790., .193634) (75990., .192734) (76190., .193309) (76390., .195187) (76590., .197916) (76790., .200849) (76990., .203277) (77190., .204611) (77390., .204545) (77590., .203145) (77790., .200820) (77990., .198179) (78190., .195861) (78390., .194379) (78590., .194028) (78790., .194843) (78990., .196608) (79190., .198895) (79390., .201158) (79590., .202850) (79790., .203568) (79990., .203165) (80190., .201783) (80390., .199804) (80590., .197730) (80790., .196050) (80990., .195132) (81190., .195151) (81390., .196070) (81590., .197653) (81790., .199517) (81990., .201214) (82190., .202338) (82390., .202629) (82590., .202044) (82790., .200761) (82990., .199120) (83190., .197528) (83390., .196352) (83590., .195846) (83790., .196098) (83990., .197026) (84190., .198396) (84390., .199876) (84590., .201110) (84790., .201808) (84990., .201815) (85190., .201155) (85390., .200014) (85590., .198686) (85790., .197493) (85990., .196707) (86190., .196492) (86390., .196876) (86590., .197752) (86790., .198901) (86990., .200046) (87190., .200912) (87390., .201300) (87590., .201131) (87790., .200467) (87990., .199486) (88190., .198436) (88390., .197567) (88590., .197073) (88790., .197055) (88990., .197499) (89190., .198288) (89390., .199226) (89590., .200088) (89790., .200669) (89990., .200840) (90190., .200572) (90390., .199946) (90590., .199127) (90790., .198316) (90990., .197704) (91190., .197425) (91390., .197531) (91590., .197985) (91790., .198670) (91990., .199417) (92190., .200047) (92390., .200414) (92590., .200437) (92790., .200125) (92990., .199563) (93190., .198896) (93390., .198285) (93590., .197872) (93790., .197745) (93990., .197924) (94190., .198356) (94390., .198932) (94590., .199511) (94790., .199957) (94990., .200166) (95190., .200097) (95390., .199777) (95590., .199290) (95790., .198760) (95990., .198312) (96190., .198049) (96390., .198026) (96590., .198239) (96790., .198630) (96990., .199102) (97190., .199540) (97390., .199842) (97590., .199940) (97790., .199817) (97990., .199512) (98190., .199103) (98390., .198691) (98590., .198374) (98790., .198222) (98990., .198264) (99190., .198485) (99390., .198826) (99590., .199203) (99790., .199526) (99990., .199719)};
\end{axis}
\end{tikzpicture}

\begin{tikzpicture}
\begin{axis}[
	xmin=0,
	xmax=100000,
	ymin=2.49,
	ymax=3.55,
	xtick={0,50000,100000,150000,200000},
	xticklabels={{\small $0$},{\small $50,000$},{\small $100,000$},{\small $150,000$},{\small $200,000$}},
	ytick={2.6,2.9,3.2,3.5},
	yticklabels={{\small $2.6\cdot 10^7$},{\small $2.9\cdot 10^7$},{\small $3.2\cdot 10^7$},{\small $3.5\cdot 10^7$}},
	axis on top=true,
	samples=300,
	xlabel=$t$,
	ylabel=$S_1(t)$,
	width=6.75cm,
	height=6.75cm,
	scaled x ticks=false,
	ylabel near ticks
]
\addplot [very thick,smooth,blue] coordinates {(190., 3.48594) (390., 3.52595) (590., 3.49970) (790., 3.47360) (990., 3.44773) (1190., 3.42207) (1390., 3.39663) (1590., 3.37140) (1790., 3.34638) (1990., 3.32158) (2190., 3.29698) (2390., 3.27259) (2590., 3.24840) (2790., 3.22442) (2990., 3.20063) (3190., 3.17705) (3390., 3.15367) (3590., 3.13048) (3790., 3.10748) (3990., 3.08468) (4190., 3.06207) (4390., 3.03965) (4590., 3.01742) (4790., 2.99538) (4990., 2.97352) (5190., 2.95184) (5390., 2.93035) (5590., 2.90903) (5790., 2.88790) (5990., 2.86694) (6190., 2.84616) (6390., 2.82555) (6590., 2.80512) (6790., 2.78486) (6990., 2.76477) (7190., 2.74484) (7390., 2.72509) (7590., 2.70550) (7790., 2.68607) (7990., 2.66681) (8190., 2.64771) (8390., 2.62877) (8590., 2.60999) (8790., 2.59136) (8990., 2.57290) (9190., 2.55458) (9390., 2.53643) (9590., 2.51843) (9790., 2.50255) (9990., 2.68308) (10190., 3.07033) (10390., 3.06356) (10590., 3.04151) (10790., 3.01928) (10990., 2.99722) (11190., 2.97535) (11390., 2.95365) (11590., 2.93215) (11790., 2.91082) (11990., 2.88967) (12190., 2.86870) (12390., 2.84790) (12590., 2.82728) (12790., 2.80683) (12990., 2.78655) (13190., 2.76645) (13390., 2.74651) (13590., 2.72674) (13790., 2.70714) (13990., 2.68770) (14190., 2.66842) (14390., 2.64931) (14590., 2.63035) (14790., 2.61156) (14990., 2.59292) (15190., 2.57444) (15390., 2.55612) (15590., 2.53797) (15790., 2.52043) (15990., 2.51261) (16190., 2.64991) (16390., 2.88162) (16590., 2.89543) (16790., 2.87724) (16990., 2.85666) (17190., 2.83600) (17390., 2.81549) (17590., 2.79514) (17790., 2.77496) (17990., 2.75495) (18190., 2.73511) (18390., 2.71544) (18590., 2.69593) (18790., 2.67659) (18990., 2.65740) (19190., 2.63838) (19390., 2.61952) (19590., 2.60082) (19790., 2.58228) (19990., 2.56393) (20190., 2.54593) (20390., 2.52953) (20590., 2.52591) (20790., 2.60732) (20990., 2.76882) (21190., 2.80619) (21390., 2.79533) (21590., 2.77686) (21790., 2.75723) (21990., 2.73749) (22190., 2.71785) (22390., 2.69834) (22590., 2.67900) (22790., 2.65981) (22990., 2.64078) (23190., 2.62193) (23390., 2.60325) (23590., 2.58479) (23790., 2.56670) (23990., 2.54959) (24190., 2.53657) (24390., 2.54355) (24590., 2.61769) (24790., 2.71948) (24990., 2.74899) (25190., 2.74182) (25390., 2.72556) (25590., 2.70715) (25790., 2.68824) (25990., 2.66925) (26190., 2.65035) (26390., 2.63160) (26590., 2.61304) (26790., 2.59478) (26990., 2.57703) (27190., 2.56049) (27390., 2.54762) (27590., 2.54736) (27790., 2.58335) (27990., 2.65743) (28190., 2.70360) (28390., 2.70964) (28590., 2.69908) (28790., 2.68327) (28990., 2.66580) (29190., 2.64781) (29390., 2.62976) (29590., 2.61187) (29790., 2.59439) (29990., 2.57776) (30190., 2.56311) (30390., 2.55359) (30590., 2.55780) (30790., 2.58992) (30990., 2.64230) (31190., 2.67596) (31390., 2.68218) (31590., 2.67410) (31790., 2.66038) (31990., 2.64451) (32190., 2.62791) (32390., 2.61123) (32590., 2.59496) (32790., 2.57976) (32990., 2.56696) (33190., 2.55962) (33390., 2.56426) (33590., 2.58870) (33790., 2.62640) (33990., 2.65371) (34190., 2.66152) (34390., 2.65647) (34590., 2.64526) (34790., 2.63140) (34990., 2.61655) (35190., 2.60162) (35390., 2.58741) (35590., 2.57501) (35790., 2.56642) (35990., 2.56541) (36190., 2.57713) (36390., 2.60203) (36590., 2.62837) (36790., 2.64309) (36990., 2.64490) (37190., 2.63847) (37390., 2.62781) (37590., 2.61525) (37790., 2.60216) (37990., 2.58953) (38190., 2.57851) (38390., 2.57081) (38590., 2.56918) (38790., 2.57690) (38990., 2.59438) (39190., 2.61484) (39390., 2.62884) (39590., 2.63300) (39790., 2.62947) (39990., 2.62135) (40190., 2.61095) (40390., 2.59976) (40590., 2.58893) (40790., 2.57967) (40990., 2.57356) (41190., 2.57275) (41390., 2.57926) (41590., 2.59269) (41790., 2.60818) (41990., 2.61933) (42190., 2.62326) (42390., 2.62092) (42590., 2.61447) (42790., 2.60582) (42990., 2.59639) (43190., 2.58737) (43390., 2.57996) (43590., 2.57561) (43790., 2.57594) (43990., 2.58203) (44190., 2.59283) (44390., 2.60452) (44590., 2.61278) (44790., 2.61564) (44990., 2.61361) (45190., 2.60821) (45390., 2.60095) (45590., 2.59310) (45790., 2.58579) (45990., 2.58016) (46190., 2.57745) (46390., 2.57881) (46590., 2.58455) (46790., 2.59330) (46990., 2.60206) (47190., 2.60791) (47390., 2.60963) (47590., 2.60759) (47790., 2.60291) (47990., 2.59679) (48190., 2.59034) (48390., 2.58459) (48590., 2.58056) (48790., 2.57924) (48990., 2.58131) (49190., 2.58659) (49390., 2.59361) (49590., 2.60007) (49790., 2.60405) (49990., 2.60481) (50190., 2.60271) (50390., 2.59860) (50590., 2.59347) (50790., 2.58826) (50990., 2.58389) (51190., 2.58122) (51390., 2.58096) (51590., 2.58342) (51790., 2.58809) (51990., 2.59362) (52190., 2.59829) (52390., 2.60086) (52590., 2.60091) (52790., 2.59880) (52990., 2.59520) (53190., 2.59094) (53390., 2.58683) (53590., 2.58365) (53790., 2.58206) (53990., 2.58254) (54190., 2.58509) (54390., 2.58908) (54590., 2.59335) (54790., 2.59666) (54990., 2.59819) (55190., 2.59776) (55390., 2.59570) (55590., 2.59259) (55790., 2.58912) (55990., 2.58596) (56190., 2.58375) (56390., 2.58299) (56590., 2.58392) (56790., 2.58636) (56990., 2.58967) (57190., 2.59289) (57390., 2.59516) (57590., 2.59595) (57790., 2.59522) (57990., 2.59329) (58190., 2.59064) (58390., 2.58786) (58590., 2.58550) (58790., 2.58407) (58990., 2.58390) (59190., 2.58506) (59390., 2.58727) (59590., 2.58993) (59790., 2.59232) (59990., 2.59380) (60190., 2.59408) (60390., 2.59320) (60590., 2.59144) (60790., 2.58923) (60990., 2.58705) (61190., 2.58536) (61390., 2.58452) (61590., 2.58473) (61790., 2.58596) (61990., 2.58788) (62190., 2.58997) (62390., 2.59169) (62590., 2.59260) (62790., 2.59254) (62990., 2.59160) (63190., 2.59005) (63390., 2.58824) (63590., 2.58658) (63790., 2.58542) (63990., 2.58501) (64190., 2.58546) (64390., 2.58664) (64590., 2.58825) (64790., 2.58986) (64990., 2.59105) (65190., 2.59155) (65390., 2.59128) (65590., 2.59036) (65790., 2.58902) (65990., 2.58758) (66190., 2.58635) (66390., 2.58560) (66590., 2.58550) (66790., 2.58606) (66990., 2.58713) (67190., 2.58845) (67390., 2.58965) (67590., 2.59044) (67790., 2.59066) (67990., 2.59027) (68190., 2.58941) (68390., 2.58829) (68590., 2.58716) (68790., 2.58628) (68990., 2.58584) (69190., 2.58593) (69390., 2.58653) (69590., 2.58747) (69790., 2.58851) (69990., 2.58939) (70190., 2.58989) (70390., 2.58991) (70590., 2.58948) (70790., 2.58870) (70990., 2.58778) (71190., 2.58693) (71390., 2.58632) (71590., 2.58610) (71790., 2.58631) (71990., 2.58689) (72190., 2.58769) (72390., 2.58849) (72590., 2.58911) (72790., 2.58940) (72990., 2.58930) (73190., 2.58885) (73390., 2.58818) (73590., 2.58745) (73790., 2.58681) (73990., 2.58642) (74190., 2.58636) (74390., 2.58663) (74590., 2.58716) (74790., 2.58781) (74990., 2.58842) (75190., 2.58884) (75390., 2.58897) (75590., 2.58880) (75790., 2.58838) (75990., 2.58781) (76190., 2.58724) (76390., 2.58678) (76590., 2.58655) (76790., 2.58658) (76990., 2.58687) (77190., 2.58734) (77390., 2.58786) (77590., 2.58831) (77790., 2.58858) (77990., 2.58861) (78190., 2.58840) (78390., 2.58802) (78590., 2.58756) (78790., 2.58712) (78990., 2.58680) (79190., 2.58668) (79390., 2.58678) (79590., 2.58707) (79790., 2.58746) (79990., 2.58787) (80190., 2.58819) (80390., 2.58834) (80590., 2.58831) (80790., 2.58809) (80990., 2.58776) (81190., 2.58739) (81390., 2.58706) (81590., 2.58686) (81790., 2.58682) (81990., 2.58694) (82190., 2.58721) (82390., 2.58753) (82590., 2.58784) (82790., 2.58806) (82990., 2.58814) (83190., 2.58806) (83390., 2.58786) (83590., 2.58757) (83790., 2.58728) (83990., 2.58705) (84190., 2.58692) (84390., 2.58693) (84590., 2.58708) (84790., 2.58731) (84990., 2.58757) (85190., 2.58780) (85390., 2.58794) (85590., 2.58796) (85790., 2.58787) (85990., 2.58768) (86190., 2.58744) (86390., 2.58722) (86590., 2.58706) (86790., 2.58699) (86990., 2.58704) (87190., 2.58718) (87390., 2.58737) (87590., 2.58758) (87790., 2.58774) (87990., 2.58783) (88190., 2.58782) (88390., 2.58771) (88590., 2.58755) (88790., 2.58736) (88990., 2.58719) (89190., 2.58709) (89390., 2.58706) (89590., 2.58712) (89790., 2.58725) (89990., 2.58741) (90190., 2.58757) (90390., 2.58768) (90590., 2.58773) (90790., 2.58770) (90990., 2.58760) (91190., 2.58745) (91390., 2.58731) (91590., 2.58719) (91790., 2.58712) (91990., 2.58712) (92190., 2.58719) (92390., 2.58730) (92590., 2.58744) (92790., 2.58755) (92990., 2.58763) (93190., 2.58764) (93390., 2.58760) (93590., 2.58751) (93790., 2.58739) (93990., 2.58728) (94190., 2.58719) (94390., 2.58716) (94590., 2.58717) (94790., 2.58724) (94990., 2.58734) (95190., 2.58744) (95390., 2.58753) (95590., 2.58757) (95790., 2.58757) (95990., 2.58752) (96190., 2.58744) (96390., 2.58735) (96590., 2.58726) (96790., 2.58720) (96990., 2.58719) (97190., 2.58722) (97390., 2.58728) (97590., 2.58736) (97790., 2.58744) (97990., 2.58750) (98190., 2.58753) (98390., 2.58751) (98590., 2.58746) (98790., 2.58739) (98990., 2.58732) (99190., 2.58726) (99390., 2.58722) (99590., 2.58722) (99790., 2.58725) (99990., 2.58731)};
\end{axis}
\end{tikzpicture}\,\,\,\begin{tikzpicture}
\begin{axis}[
	xmin=0,
	xmax=100000,
	ymin=0,
	ymax=8.9,
	xtick={0,50000,100000,150000,200000},
	xticklabels={{\small $0$},{\small $50,000$},{\small $100,000$},{\small $150,000$},{\small $200,000$}},
	ytick={2,4,6,8},
	yticklabels={{\small $2\cdot 10^6$},{\small $4\cdot 10^6$},{\small $6\cdot 10^6$},{\small $8\cdot 10^6$}},
	axis on top=true,
	samples=300,
	xlabel=$t$,
	ylabel=$I_2(t)$,
	width=6.75cm,
	height=6.75cm,
	scaled x ticks=false,
	ylabel near ticks
]
\addplot [very thick,smooth,red] coordinates {(190., 8.88826) (390., 0.128599e-1) (590., 0) (790., 0) (990., 0) (1190., 0) (1390., 0) (1590., 0) (1790., 0) (1990., 0) (2190., 0) (2390., 0) (2590., 0) (2790., 0) (2990., 0) (3190., 0) (3390., 0) (3590., 0) (3790., 0) (3990., 0) (4190., 0) (4390., 0) (4590., 0) (4790., 0) (4990., 0) (5190., 0) (5390., 0) (5590., 0) (5790., 0) (5990., 0) (6190., 0) (6390., 0) (6590., 0) (6790., 0) (6990., 0) (7190., 0) (7390., 0) (7590., 0) (7790., 0) (7990., 0) (8190., 0) (8390., 0) (8590., 0) (8790., 0) (8990., 0) (9190., 0) (9390., 0) (9590., 0) (9790., 0.139909e-1) (9990., 1.33124) (10190., .411104) (10390., 0) (10590., 0) (10790., 0) (10990., 0) (11190., 0) (11390., 0) (11590., 0) (11790., 0) (11990., 0) (12190., 0) (12390., 0) (12590., 0) (12790., 0) (12990., 0) (13190., 0) (13390., 0) (13590., 0) (13790., 0) (13990., 0) (14190., 0) (14390., 0) (14590., 0) (14790., 0) (14990., 0) (15190., 0) (15390., 0) (15590., 0) (15790., 0) (15990., 0.535659e-1) (16190., .724654) (16390., .400841) (16590., 0.342362e-1) (16790., 0) (16990., 0) (17190., 0) (17390., 0) (17590., 0) (17790., 0) (17990., 0) (18190., 0) (18390., 0) (18590., 0) (18790., 0) (18990., 0) (19190., 0) (19390., 0) (19590., 0) (19790., 0) (19990., 0) (20190., 0) (20390., 0) (20590., 0.630764e-1) (20790., .390727) (20990., .383291) (21190., 0.759859e-1) (21390., 0.124378e-1) (21590., 0) (21790., 0) (21990., 0) (22190., 0) (22390., 0) (22590., 0) (22790., 0) (22990., 0) (23190., 0) (23390., 0) (23590., 0) (23790., 0) (23990., 0) (24190., 0.194391e-1) (24390., 0.953029e-1) (24590., .302602) (24790., .254398) (24990., 0.765241e-1) (25190., 0.194066e-1) (25390., 0) (25590., 0) (25790., 0) (25990., 0) (26190., 0) (26390., 0) (26590., 0) (26790., 0) (26990., 0) (27190., 0) (27390., 0.182652e-1) (27590., 0.618366e-1) (27790., .174318) (27990., .231813) (28190., .121921) (28390., 0.430725e-1) (28590., 0.155461e-1) (28790., 0) (28990., 0) (29190., 0) (29390., 0) (29590., 0) (29790., 0) (29990., 0) (30190., 0.113786e-1) (30390., 0.277749e-1) (30590., 0.712447e-1) (30790., .148198) (30990., .171664) (31190., .104004) (31390., 0.464912e-1) (31590., 0.208212e-1) (31790., 0.107657e-1) (31990., 0) (32190., 0) (32390., 0) (32590., 0) (32790., 0) (32990., 0.162021e-1) (33190., 0.327815e-1) (33390., 0.683066e-1) (33590., .119476) (33790., .136299) (33990., 0.962507e-1) (34190., 0.521095e-1) (34390., 0.273639e-1) (34590., 0.159562e-1) (34790., 0.109816e-1) (34990., 0) (35190., 0) (35390., 0.113219e-1) (35590., 0.164902e-1) (35790., 0.275865e-1) (35990., 0.494353e-1) (36190., 0.834624e-1) (36390., .110795) (36590., .102257) (36790., 0.697004e-1) (36990., 0.418100e-1) (37190., 0.256724e-1) (37390., 0.176131e-1) (37590., 0.140812e-1) (37790., 0.133303e-1) (37990., 0.149428e-1) (38190., 0.195649e-1) (38390., 0.290038e-1) (38590., 0.459112e-1) (38790., 0.703138e-1) (38990., 0.913684e-1) (39190., 0.911537e-1) (39390., 0.708180e-1) (39590., 0.479815e-1) (39790., 0.321033e-1) (39990., 0.231230e-1) (40190., 0.188024e-1) (40390., 0.176256e-1) (40590., 0.191018e-1) (40790., 0.236380e-1) (40990., 0.324176e-1) (41190., 0.467039e-1) (41390., 0.651167e-1) (41590., 0.796247e-1) (41790., 0.796696e-1) (41990., 0.658213e-1) (42190., 0.485877e-1) (42390., 0.352259e-1) (42590., 0.269737e-1) (42790., 0.227920e-1) (42990., 0.216918e-1) (43190., 0.233192e-1) (43390., 0.279576e-1) (43590., 0.362882e-1) (43790., 0.485207e-1) (43990., 0.623755e-1) (44190., 0.716593e-1) (44390., 0.704269e-1) (44590., 0.599969e-1) (44790., 0.470390e-1) (44990., 0.364600e-1) (45190., 0.295938e-1) (45390., 0.260846e-1) (45590., 0.253967e-1) (45790., 0.273388e-1) (45990., 0.320863e-1) (46190., 0.398649e-1) (46390., 0.501019e-1) (46590., 0.601832e-1) (46790., 0.655724e-1) (46990., 0.632494e-1) (47190., 0.549465e-1) (47390., 0.450900e-1) (47590., 0.369139e-1) (47790., 0.315176e-1) (47990., 0.288717e-1) (48190., 0.287335e-1) (48390., 0.310121e-1) (48590., 0.357477e-1) (48790., 0.427662e-1) (48990., 0.510065e-1) (49190., 0.579869e-1) (49390., 0.606432e-1) (49590., 0.577234e-1) (49790., 0.509740e-1) (49990., 0.434054e-1) (50190., 0.371654e-1) (50390., 0.330969e-1) (50590., 0.313081e-1) (50790., 0.317119e-1) (50990., 0.342426e-1) (51190., 0.387849e-1) (51390., 0.448591e-1) (51590., 0.512016e-1) (51790., 0.557453e-1) (51990., 0.566037e-1) (52190., 0.535063e-1) (52390., 0.479986e-1) (52590., 0.421630e-1) (52790., 0.374547e-1) (52990., 0.345093e-1) (53190., 0.334762e-1) (53390., 0.343292e-1) (53590., 0.369727e-1) (53790., 0.411442e-1) (53990., 0.461796e-1) (54190., 0.508387e-1) (54390., 0.535662e-1) (54590., 0.533231e-1) (54790., 0.503243e-1) (54990., 0.458466e-1) (55190., 0.413530e-1) (55390., 0.378528e-1) (55590., 0.358234e-1) (55790., 0.354062e-1) (55990., 0.365775e-1) (56190., 0.391850e-1) (56390., 0.428553e-1) (56590., 0.468585e-1) (56790., 0.501176e-1) (56990., 0.515593e-1) (57190., 0.507004e-1) (57390., 0.479558e-1) (57590., 0.443423e-1) (57790., 0.409038e-1) (57990., 0.383554e-1) (58190., 0.370523e-1) (58390., 0.371028e-1) (58590., 0.384556e-1) (58790., 0.409002e-1) (58990., 0.440021e-1) (59190., 0.470596e-1) (59390., 0.492169e-1) (59590., 0.497924e-1) (59790., 0.486392e-1) (59990., 0.462213e-1) (60190., 0.433335e-1) (60390., 0.407317e-1) (60590., 0.389295e-1) (60790., 0.381856e-1) (60990., 0.385643e-1) (61190., 0.399763e-1) (61390., 0.421686e-1) (61590., 0.446947e-1) (61790., 0.469390e-1) (61990., 0.482691e-1) (62190., 0.482918e-1) (62390., 0.470470e-1) (62590., 0.449764e-1) (62790., 0.426959e-1) (62990., 0.407589e-1) (63190., 0.395361e-1) (63390., 0.392075e-1) (63590., 0.397926e-1) (63790., 0.411674e-1) (63990., 0.430583e-1) (64190., 0.450453e-1) (64390., 0.466252e-1) (64590., 0.473599e-1) (64790., 0.470535e-1) (64990., 0.458388e-1) (65190., 0.441059e-1) (65390., 0.423305e-1) (65590., 0.409198e-1) (65790., 0.401397e-1) (65990., 0.401058e-1) (66190., 0.407980e-1) (66390., 0.420680e-1) (66590., 0.436433e-1) (66790., 0.451548e-1) (66990., 0.462143e-1) (67190., 0.465379e-1) (67390., 0.460563e-1) (67590., 0.449396e-1) (67790., 0.435185e-1) (67990., 0.421593e-1) (68190., 0.411618e-1) (68390., 0.407121e-1) (68590., 0.408748e-1) (68790., 0.415984e-1) (68990., 0.427229e-1) (69190., 0.439945e-1) (69390., 0.451058e-1) (69590., 0.457726e-1) (69790., 0.458252e-1) (69990., 0.452707e-1) (70190., 0.442851e-1) (70390., 0.431421e-1) (70590., 0.421222e-1) (70790., 0.414444e-1) (70990., 0.412341e-1) (71190., 0.415159e-1) (71390., 0.422173e-1) (71590., 0.431782e-1) (71790., 0.441745e-1) (71990., 0.449619e-1) (72190., 0.453424e-1) (72390., 0.452270e-1) (72590., 0.446648e-1) (72790., 0.438217e-1) (72990., 0.429200e-1) (73190., 0.421730e-1) (73390., 0.417382e-1) (73590., 0.416944e-1) (73790., 0.420368e-1) (73990., 0.426813e-1) (74190., 0.434772e-1) (74390., 0.442352e-1) (74590., 0.447691e-1) (74790., 0.449484e-1) (74990., 0.447381e-1) (75190., 0.442077e-1) (75390., 0.435050e-1) (75590., 0.428080e-1) (75790., 0.422771e-1) (75990., 0.420225e-1) (76190., 0.420882e-1) (76390., 0.424491e-1) (76590., 0.430174e-1) (76790., 0.436582e-1) (76990., 0.442174e-1) (77190., 0.445589e-1) (77390., 0.446028e-1) (77590., 0.443480e-1) (77790., 0.438710e-1) (77990., 0.432989e-1) (78190., 0.427720e-1) (78390., 0.424092e-1) (78590., 0.422844e-1) (78790., 0.424160e-1) (78990., 0.427667e-1) (79190., 0.432511e-1) (79390., 0.437532e-1) (79590., 0.441520e-1) (79790., 0.443520e-1) (79990., 0.443095e-1) (80190., 0.440439e-1) (80390., 0.436300e-1) (80590., 0.431744e-1) (80790., 0.427862e-1) (80990., 0.425512e-1) (81190., 0.425163e-1) (81390., 0.426820e-1) (81590., 0.430045e-1) (81790., 0.434055e-1) (81990., 0.437883e-1) (82190., 0.440613e-1) (82390., 0.441606e-1) (82590., 0.440675e-1) (82790., 0.438122e-1) (82990., 0.434634e-1) (83190., 0.431087e-1) (83390., 0.428311e-1) (83590., 0.426911e-1) (83790., 0.427153e-1) (83990., 0.428923e-1) (84190., 0.431770e-1) (84390., 0.435001e-1) (84590., 0.437837e-1) (84790., 0.439608e-1) (84990., 0.439912e-1) (85190., 0.438728e-1) (85390., 0.436400e-1) (85590., 0.433536e-1) (85790., 0.430837e-1) (85990., 0.428926e-1) (86190., 0.428212e-1) (86390., 0.428812e-1) (86590., 0.430545e-1) (86790., 0.432975e-1) (86990., 0.435511e-1) (87190., 0.437547e-1) (87390., 0.438606e-1) (87590., 0.438464e-1) (87790., 0.437198e-1) (87990., 0.435157e-1) (88190., 0.432860e-1) (88390., 0.430860e-1) (88590., 0.429610e-1) (88790., 0.429370e-1) (88990., 0.430160e-1) (89190., 0.431762e-1) (89390., 0.433777e-1) (89590., 0.435717e-1) (89790., 0.437123e-1) (89990., 0.437671e-1) (90190., 0.437260e-1) (90390., 0.436023e-1) (90590., 0.434289e-1) (90790., 0.432490e-1) (90990., 0.431051e-1) (91190., 0.430294e-1) (91390., 0.430368e-1) (91590., 0.431228e-1) (91790., 0.432647e-1) (91990., 0.434275e-1) (92190., 0.435719e-1) (92390., 0.436642e-1) (92590., 0.436838e-1) (92790., 0.436285e-1) (92990., 0.435145e-1) (93190., 0.433711e-1) (93390., 0.432335e-1) (93590., 0.431337e-1) (93790., 0.430936e-1) (93990., 0.431203e-1) (94190., 0.432054e-1) (94390., 0.433269e-1) (94590., 0.434551e-1) (94790., 0.435594e-1) (94990., 0.436156e-1) (95190., 0.436120e-1) (95390., 0.435514e-1) (95590., 0.434506e-1) (95790., 0.433349e-1) (95990., 0.432323e-1) (96190., 0.431664e-1) (96390., 0.431512e-1) (96590., 0.431884e-1) (96790., 0.432676e-1) (96990., 0.433687e-1) (97190., 0.434671e-1) (97390., 0.435397e-1) (97590., 0.435699e-1) (97790., 0.435519e-1) (97990., 0.434919e-1) (98190., 0.434056e-1) (98390., 0.433145e-1) (98590., 0.432403e-1) (98790., 0.431997e-1) (98990., 0.432010e-1) (99190., 0.432425e-1) (99390., 0.433130e-1) (99590., 0.433950e-1) (99790., 0.434686e-1) (99990., 0.435167e-1)};
\end{axis}
\end{tikzpicture}

\caption{\label{fig:TShigh} Plots of $S(t)$, $I_1(t)$, $S_1(t)$, and $I_2(t)$, versus $t$, for $0\leqslant t\leqslant 100,000$, in the high-transmission case.}
\end{figure}

For the parameter values presented in Table \ref{tab:parameters}, $\upsilon=0.0167$, and $\beta=8\cdot 10^{-10}$, we have that $\cR_0\approx 4.62570$ and
$$\mathbf{e}^{(1)}=\left(S^{(1)},I_1^{(1)},S_1^{(1)}, I_2^{(1)}\right)\approx\left(4.12\cdot 10^7,\, 1.99\cdot 10^5,\, 2.59\cdot 10^7,\, 4.34\cdot 10^{4}\right).$$
Figure \ref{fig:TShigh} shows the convergence of the model's solution associated to the initial condition given by \eqref{eq:IC} to the endemic equilibrium $\mathbf{e}^{(1)}$. Notice not only the presence of oscillations as in the low-transmission case, but also the significantly faster convergence rate, the shorter period of oscillations, and the considerably higher values achieved by $I_2(t)$ compared to those in the low-transmission case.

\begin{table}\renewcommand{\arraystretch}{1.2}
\begin{tabular}{|c|r|r|r|r|}\hline
$p$ & \multicolumn{1}{c|}{$\Upsilon _{p}^{S^{(1)}}$} & \multicolumn{1}{c|}{$\Upsilon _{p}^{I_1^{(1)}}$} & \multicolumn{1}{c|}{$\Upsilon _{p}^{S_1^{(1)}}$} & \multicolumn{1}{c|}{$\Upsilon _{p}^{I_2^{(1)}}$}\\ \hhline{|=|=|=|=|=|}
$\upsilon $ & $-0.01231$ & $-0.01800$ & $0.05648$ & $0.05078$ \\\hline
$\mu $ & $0.11024$ & $-0.24237$ & $-0.50004$ & $-0.85265$ \\\hline
$\beta $ & $-1.10884$ & $0.24102$ & $-0.50052$ & $0.84934$ \\\hline
$\mu'$ & $0.96670$ & $-0.93504$ & $-0.38459$ & $-2.28634$ \\\hline
$\alpha $ & $0.14074$ & $-0.30462$ & $0.88567$ & $0.44032$ \\\hline
$\delta $ & $-0.09190$ & $0.01997$ & $-0.57830$ & $0.53357$ \\\hline
$\Lambda $  & $-0.10884$ & $1.24102$ & $0.49948$ & $1.84934$\\\hline
\end{tabular}\smallskip
\caption{\label{tab:sensitivityhigh}The sensitivity indices of $S^{(1)}$, $I_1^{(1)}$, $S_1^{(1)}$, and $I_2^{(1)}$ with respect to $p\in\left\{\upsilon,\Lambda,\mu,\beta,\mu',\alpha,\delta\right\}$ in the high-transmission case.}
\end{table}

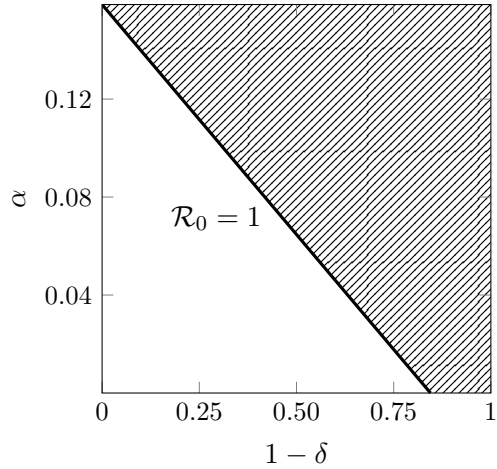
\begin{figure}\centering
\begin{tikzpicture}
\pgfset{declare function={f(\x)=9406408899/59312500000-11733301/62500000*\x;}}

\begin{axis}[
	xmin=0,
	xmax=1,
	ymin=0,
	ymax=0.1586,
	xtick={0,0.25,0.50,0.75,1},
	xticklabels={{\small $0$},{\small $0.25$},{\small $0.50$},{\small $0.75$},{\small $1$}},
	ytick={0.04,0.08,0.12},
	yticklabels={{\small $0.04$},{\small $0.08$},{\small $0.12$}},
	axis on top=true,
	samples=100,
	xlabel=$1-\delta$,
	ylabel=$\alpha$,
	width=6.75cm,
	height=6.75cm,
	scaled x ticks=false,
	ylabel near ticks
]
\addplot [very thick,domain=0:0.84477] {f(x)};
\fill[pattern=north east lines] (axis cs:0,0.1586) -- (axis cs:0.84477,0) -- (axis cs:1,0) -- (axis cs:1,0.1586) -- cycle;

\node[below left,xshift=2pt,yshift=1pt] at (axis cs:.422385,0.79295e-1) {$\mathcal{R}_0=1$};


\end{axis}
\end{tikzpicture}
\caption{\label{fig:R0=1extreme}The curve $\cR_0=1$ on the $(1-\delta)\alpha$-plane, with the disease-free region shaded.}
\end{figure}

\begin{figure}\centering
	\includegraphics[width=0.5\textwidth]{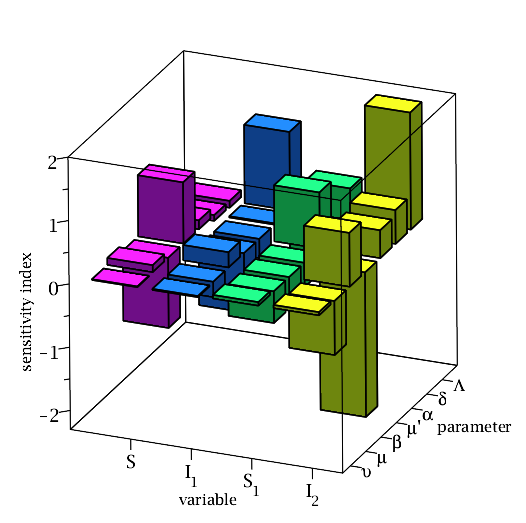}
	\caption{\label{fig:sensitivityhigh}Three-dimensional boxplot of the sensitivity indices of $S^{(1)}$, $I_1^{(1)}$, $S_1^{(1)}$, and $I_2^{(1)}$ with respect to $p\in\left\{\upsilon,\Lambda,\mu,\beta,\mu',\alpha,\delta\right\}$ in the high-transmission case.}
\end{figure}

In this case, the values of the sensitivity indices are presented in Table \ref{tab:sensitivityhigh} and visualised in Figure \ref{fig:sensitivityhigh}. Notice, in particular, the drastically lower values of $\Upsilon_p^{I_1^{(1)}}$ and $\Upsilon_p^{I_2^{(1)}}$ for $p\in\left\{\Lambda,\mu,\beta\right\}$ compared to those in the low-transmission case, which provide a justification for the difficulty of suppressing the value of $\cR_0$ in the situation where transmissions take place at such a high rate. Indeed, from \eqref{eq:R0dlmbetaupsilon} one readily finds that, in the case of $\upsilon=100\%$, $\cR_0<1$ corresponds to $\beta<4.93\cdot 10^{-10}$. This means that, with such a high transmission rate, even achieving four-dose vaccination across the entire population is insufficient. The effort of maximising the proportion of vaccinated individuals, therefore, should be accompanied by other initiatives, namely, improvements in the vaccine's efficacy $1-\delta$ and the disease's recovery rate $\alpha$. In the extreme case $\upsilon=1$, from \eqref{eq:R0} have that
\begin{equation}\label{eq:R0dlm1-deltaalpha}
\cR_0\left(1-\delta,\alpha\right)=\frac{11134902649\left[1-\left(1-\delta\right)\right]}{1728493750+59312500000\alpha}.
\end{equation}
The curve $\cR_0\left(1-\delta,\alpha\right)=1$ is plotted in Figure \ref{fig:R0=1extreme}, where, as before, the shaded region is the disease-free region, where $\cR_0\left(1-\delta,\alpha\right)<1$. From \eqref{eq:R0dlm1-deltaalpha} one finds that, without changing the present recovery rate $\alpha=0.011$, the new variants' transmission risks can be mitigated only if the vaccine's efficacy $1-\delta$ is raised at least to $78.62\%$. On the other hand, without changing the present vaccine efficacy $1-\delta=65.3\%$, the same risks can be mitigated by increasing the recovery rate at least to $0.036$ per day.

\section{Conclusions and future research}\label{sec:conclusions}

We have constructed a mathematical model as a means to assess the effectiveness of four-dose COVID-19 vaccinations in countering transmissions of a novel virus variant. We have computed the model's basic reproduction number $\mathcal{R}_0$, and shown that the model has a unique disease-free equilibrium, which exists for all parameter values and is stable if $\cR_0<1$. In addition, we have proved that a unique positive endemic equilibrium exists if and only if $\mathcal{R}_0>1$.

Using a set of parameter values specified to represent the situation in Indonesia as of March 20, 2023, we have simulated two different cases corresponding to two different transmission levels of the new variant. Our findings indicate that, if social restrictions are not to be reinforced, an intense acceleration of four-dose vaccinations is necessary to mitigate the new variant's transmissions, with a significant risk of ineffectiveness in the case of a high transmission level. Therefore, to answer the question posed in our title, accelerating four-dose vaccinations alone is not sufficient to overcome risks of new variant transmissions. We recommend that four-dose vaccinations are not only accelerated but also operated alongside other forms of intervention, such as raising the vaccine's efficacy and the disease's recovery rate.


As a continuation of the present research, one could construct mathematical models for the spread of COVID-19 which relate the emergence of new variants to certain vaccination strategies. The fact that our model's solutions exhibit oscillations calls to mind the so-called pulse vaccination strategy \cite{AgurCojocaruMazorAndersonDanon,dOnofrio,EarnRohaniGrenfell}, which involves a periodic administration of vaccinations to specific subpopulations, and proves to be promising for diseases exhibiting periodic outbreaks \cite{dOnofrio,EarnRohaniGrenfell}. Investigations in this direction could aim to construct a sequence of time-intervals along which vaccinations must be administered for an optimal retainment of the population's immunity against COVID-19, in response to the increasingly diverse range of variants.




\section*{Disclosure statement}

\noindent No potential conflict of interest was reported by the authors.


\end{document}